\documentclass[12pt]{amsart}
\usepackage[utf8]{inputenc}
\usepackage{amsmath, amssymb, amsfonts, tikz, comment, url, wasysym, import, geometry, amsthm, mathtools, subcaption, xcolor}
\usepackage[capitalise]{cleveref}
\usepackage[export]{adjustbox}
\usepackage{cite}
\usepackage{setspace}
\usepackage{indentfirst}

\def\Z{\mathbb{Z}}

\newtheorem{theorem}{Theorem}[section]
\newtheorem{lemma}[theorem]{Lemma}
\newtheorem{proposition}[theorem]{Proposition}

\newtheorem{conjecture}[theorem]{Conjecture}
\newtheorem{question}[theorem]{Question}

\theoremstyle{definition}
\newtheorem{example}{Example}[section]
\newtheorem{definition}{Definition}[section]

\theoremstyle{remark}

\makeatletter
\renewcommand*\env@matrix[1][\arraystretch]{%
  \edef\arraystretch{#1}%
  \hskip -\arraycolsep
  \let\@ifnextchar\new@ifnextchar
  \array{*\c@MaxMatrixCols c}}
\makeatother

\def\cp{\,\begin{tikzpicture}
    \draw[black] (0,0) -- (0.25,0) -- (0.25,0.25) -- (0,0.25) -- (0,0);
    \end{tikzpicture}\,}
\def\sp{\,\begin{tikzpicture}
    \draw[black] (0,0) -- (0.25,0) -- (0.25,0.25) -- (0,0.25) -- (0,0) -- (0.25,0.25) -- (0, 0.25) -- (0.25,0);
    (\end{tikzpicture}\,}

\begin{document}

\title[Compatible Recurrent Identities of the Sandpile Group]{Compatible Recurrent Identities of the Sandpile Group and Maximal Stable Configurations}
\author{Yibo Gao}
\address{Department of Mathematics, Massachusetts Institute of Technology, \mbox{Cambridge, MA 02139}}
\email{gaoyibo@mit.edu}
\author{Rupert Li}
\address{Massachusetts Institute of Technology,
\mbox{Cambridge, MA 02139}}
\email{rupertli@mit.edu}
\keywords{chip-firing, recurrent configuration, sandpile group, reduced Laplacian, maximal stable configuration, recurrent identity, Cartesian product, strong product, bipartite graph, tree}


\date{June 23, 2020}

\begin{abstract}
    In the abelian sandpile model, recurrent chip configurations are of interest as they are a natural choice of coset representatives under the quotient of the reduced Laplacian. We investigate graphs whose recurrent identities with respect to different sinks are compatible with each other. The maximal stable configuration is the simplest recurrent chip configuration, and graphs whose recurrent identities equal the maximal stable configuration are of particular interest, and are said to have the complete maximal identity property. We prove that given any graph $G$ one can attach trees to the vertices of $G$ to yield a graph with the complete maximal identity property. We conclude with several intriguing conjectures about the complete maximal identity property of various graph products.
\end{abstract}



\maketitle

\section{Introduction}\label{section:Introduction}
    The \textbf{abelian sandpile model} was invented by Dhar in \cite{dhar1990self} as a model for self-organized criticality, which was introduced in \cite{bak1987self}.
    The automaton model was invented to model stacked items at sites, which would topple if a critical height was reached, and send the items to adjacent sites; for example, using a lattice graph to model a plane, the abelian sandpile model could mimic a pile of sand on a flat surface as it collapses under gravity to reach a certain stable state, and the destabilizing behavior of adding additional sand analyzed.
    The analysis of these models have led to the generalization of it for arbitrary locations and compositions of sites.
    
    Among the models that display self-organized criticality, the abelian sandpile model, which has been used to model landslides, is the simplest analytically tractable model, according to \cite{dhar1999abelian}.
    Overviews of the abelian sandpile model, related models, and self-organized criticality are given in \cite{dhar2006theoretical,bak2013nature}.
    It has been demonstrated in \cite{bak1990forest,burridge1967model,turcotte2004landslides} that self-organized criticality is present in models that can be associated to the natural hazards of landslides, earthquakes, and forest fires, as well as to financial markets in \cite{bartolozzi2006scale,biondo2015modeling,scheinkman1994self}, evolution and species extinction in \cite{sneppen1995evolution,newman1996self,sole1996extinction}, and neural systems in the brain in \cite{hesse2014self,pu2013developing,chialvo2010emergent}.

    The abelian sandpile model, also referred to as the \textbf{chip-firing game} after the paper \cite{bjorner1991chip}, on a directed graph $G$ consists of a collection of chips at each vertex of $G$.
    If a vertex $v$ has at least as many chips as its outdegree, then it can fire, sending one chip along each outgoing edge to its neighboring vertices.
    This continues indefinitely or until no vertex can fire.
    Conditions for which configurations lead to terminating processes have been studied in \cite{bjorner1991chip,bjorner1992chip}, and a polynomial bound for the length of the process was proven in \cite{tardos1988polynomial}.
    Chip-firing has been studied using an algebraic potential theory approach in \cite{biggs1997algebraic}; a particular type of chip-firing, referred to as the probabilistic abacus, has also been considered in \cite{engel1975probabilistic,engel1976does} as a quasirandom process that provides insight into Markov chains.
    Books surveying chip-firing include \cite{klivans2018mathematics,corry2018divisors}.
    
    In \cref{section:Preliminaries}, we establish the basic theory surrounding the chip-firing game, and define the primary algebraic object associated with the model, the \textbf{sandpile group}, and the special role \textbf{recurrent} elements play in the group.
    The sandpile group has been studied in \cite{creutz1991abelian}, as well as in \cite{cori2000sandpile} and \cite{levine2009sandpile} for the cases of dual graphs and wired trees, respectively.
    The algebraic structure of the sandpile group has been studied in \cite{wagner2000critical,dhar1995algebraic}, and the properties of recurrent configurations have been studied in \cite{biggs1999chip,van2001algorithmic}.
    
    In \cref{section:CMIP}, we begin stating our own results, investigating the graphs for which the simplest recurrent element, the \textbf{maximal stable configuration}, is the identity of the sandpile group, regardless of the choice of sink.
    The identity of the sandpile group is of particular interest, and has been previously investigated in \cite{caracciolo2008explicit,le2002identity} for $\mathbb{Z}^2$ lattice graphs.
    In \cref{section: necessary and sufficient}, we find necessary conditions for the property developed in \cref{section:CMIP}, the \textbf{complete maximal identity property}.
    In \cref{section:CIP}, we generalize the complete maximal identity property into the \textbf{complete identity property}, investigating graphs for which the recurrent identity remains the same regardless of the choice of sink.
    We consider the possible relationship between the complete identity property and bipartite graphs.
    We conclude in \cref{section:Conjectures} with some conjectures surrounding the recurrent identity of sandpiles formed from graph products. 

\section{Preliminaries}\label{section:Preliminaries}
    We refer readers to \cite{corry2018divisors} and \cite{holroyd2008chip} for detailed background on chip-firing.
    Unless otherwise stated, we restrict the chip-firing game to be on graphs which are simple, connected, and nontrivial (have more than one vertex).
    Recall that a simple graph is an unweighted, undirected graph without loops or multiple edges.
    A \textbf{sandpile} is a graph $G$ that has a special vertex, called a \textbf{sink}.
    A \textbf{chip configuration} over the sandpile is a vector of integers indexed over all non-sink vertices of $G$.
    In standard convention, these numbers must be nonnegative integers, a discrete number of \textbf{chips}.
    A chip configuration $c$, if explicitly stated to be over all vertices, not just all non-sink vertices, may be used for convenience when the sink changes.
    The entry corresponding to the sink is simply excluded from the vector.
    
    In a sandpile, a vertex can \textbf{fire} if it has at least as many chips as its degree, at which point it sends chips along each edge to its neighboring vertices, with each edge transferring $w$, where $w$ is the weight of the edge (in an undirected graph, each edge has weight 1); the vertex that fires loses the chips it fired.
    A vertex is said to be \textbf{active} if it can fire.
    The sink is not allowed to fire, nor are its chips considered; hence a chip configuration for a sandpile does not have an entry corresponding to the number of chips at the sink.
    A \textbf{stable} configuration is a configuration that has no active vertices.
    If by a sequence of firings a chip configuration $c$ can result in a stable configuration, that stable configuration is called the \textbf{stabilization} and is denoted $\operatorname{Stab}(c)$.
    It has been proven in \cite{bjorner1991chip} that if a stabilization exists, then it is unique for each chip configuration, regardless of the order in which vertices are fired; moreover, regardless of the order in which vertices are fired, the chip configuration will eventually result in the stable configuration, always taking the same number of firings in total.
    It has also been shown in \cite{holroyd2008chip} that in a sandpile all chip configurations stabilize.
    
    \begin{definition}\label{definition:MSC}
        The \textbf{maximal stable configuration} $m_G$ is the chip configuration in which every vertex $v$ has $d_v-1$ chips, where $d_v$ is the degree of vertex $v$.
    \end{definition}

    For a sandpile, a chip configuration $c$ is called \textbf{accessible} if for all (stable) configurations $d$ there exists a configuration $e$ such that $\operatorname{Stab}(d+e)=c$, where $d+e$ is the componentwise addition of the vectors $d$ and $e$.
    
    \begin{definition}\label{definition: recurrent}
    A chip configuration is called \textbf{recurrent} if it is accessible and stable.
    \end{definition}
    
    Note that $m_G$ is always recurrent, even when there is only one recurrent configuration for a given sandpile, and is by far the simplest recurrent chip configuration.
    
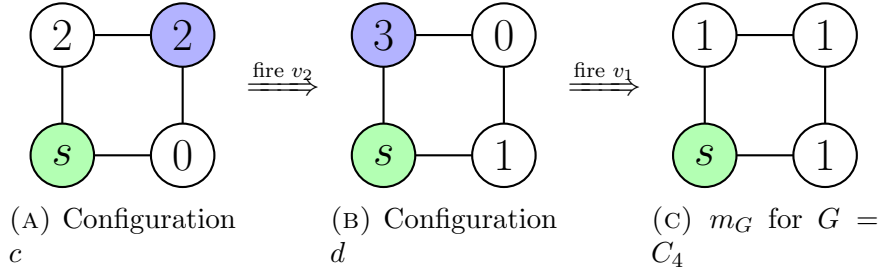
\begin{figure}[h!]
        \centering
        \begin{subfigure}[t]{0.2\textwidth}
            \centering
            \begin{tikzpicture}[scale=0.8,font=\Large,baseline]
                \draw[thick] (-1,-1) -- (1,-1) -- (1,1) -- (-1,1) -- (-1,-1);
                \filldraw[color=black,fill=green!30,thick] (-1,-1) circle (15pt);
                \filldraw[color=black,fill=white,thick] (-1,1) circle (15pt);
                \filldraw[color=black,fill=blue!30,thick] (1,1) circle (15pt);
                \filldraw[color=black,fill=white,thick] (1,-1) circle (15pt);
                \node at (-1,-1) {$s$};
                \node at (-1,1) {2};
                \node at (1,1) {2};
                \node at (1,-1) {0};
            \end{tikzpicture}
            \caption{Configuration $c$}
            \label{subfigure: C4 220}
        \end{subfigure}
        $\xRightarrow{\text{fire } v_2}$
        \begin{subfigure}[t]{0.2\textwidth}
            \centering
            \begin{tikzpicture}[scale=0.8,font=\Large,baseline]
                \draw[thick] (-1,-1) -- (1,-1) -- (1,1) -- (-1,1) -- (-1,-1);
                \filldraw[color=black,fill=green!30,thick] (-1,-1) circle (15pt);
                \filldraw[color=black,fill=blue!30,thick] (-1,1) circle (15pt);
                \filldraw[color=black,fill=white,thick] (1,1) circle (15pt);
                \filldraw[color=black,fill=white,thick] (1,-1) circle (15pt);
                \node at (-1,-1) {$s$};
                \node at (-1,1) {3};
                \node at (1,1) {0};
                \node at (1,-1) {1};
            \end{tikzpicture}
            \caption{Configuration $d$}
            \label{subfigure: C4 301}
        \end{subfigure}
        $\xRightarrow{\text{fire } v_1}$
        \begin{subfigure}[t]{0.2\textwidth}
            \centering
            \begin{tikzpicture}[scale=0.8,font=\Large,baseline]
                \draw[thick] (-1,-1) -- (1,-1) -- (1,1) -- (-1,1) -- (-1,-1);
                \filldraw[color=black,fill=green!30,thick] (-1,-1) circle (15pt);
                \filldraw[color=black,fill=white,thick] (-1,1) circle (15pt);
                \filldraw[color=black,fill=white,thick] (1,1) circle (15pt);
                \filldraw[color=black,fill=white,thick] (1,-1) circle (15pt);
                \node at (-1,-1) {$s$};
                \node at (-1,1) {1};
                \node at (1,1) {1};
                \node at (1,-1) {1};
            \end{tikzpicture}
            \caption{$m_G$ for $G = C_4$}
            \label{subfigure: C4 111}
        \end{subfigure}
        \caption{Chip-firing on $C_4$. The non-sink vertices are labeled with their number of chips.}
        \label{figure: Chip-firing example on C4}
    \end{figure}

    \begin{example}\label{example: preliminary defs C4}
        The cycle graph on four vertices $C_4$ may have its vertices indexed 0, 1, 2, and 3, where $v_0$ is the sink, as shown in \cref{subfigure: C4 Sandpile}, where the sink is colored green.
        
        \begin{figure}[htbp!]
            \centering
            \begin{tikzpicture}[scale=0.8,font=\Large,baseline]
                \draw[thick] (-1,-1) -- (1,-1) -- (1,1) -- (-1,1) -- (-1,-1);
                \filldraw[color=black,fill=white,thick] (-1,1) circle (15pt);
                \filldraw[color=black,fill=white,thick] (1,-1) circle (15pt);
                \filldraw[color=black,fill=white,thick] (1,1) circle (15pt);
                \filldraw[color=black,fill=green!30,thick] (-1,-1) circle (15pt);
                \node at (-1,-1) {$s$};
                \node at (1,-1) {$v_3$};
                \node at (1,1) {$v_2$};
                \node at (-1,1) {$v_1$};
            \end{tikzpicture}
            \caption{$C_4$ with sink $v_0$, labelled as $s$.}
            \label{subfigure: C4 Sandpile}
        \end{figure}
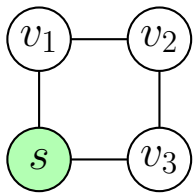
        
        This sandpile may have a chip configuration $c=(2,2,0)$ as shown in \cref{subfigure: C4 220}.
        In $c$, vertices $v_1$ and $v_2$ are active.
        If $v_2$ (colored blue in \cref{subfigure: C4 220}) fires, the result is chip configuration $d=(3,0,1)$ as shown in \cref{subfigure: C4 301}.
        In $d$, $v_1$ is the only active vertex (colored blue in \cref{subfigure: C4 301}), and firing it results in the maximal stable configuration $m_G=(1,1,1)$ where in this example $G = C_4$ with sink at $v_0$, as shown in \cref{subfigure: C4 111}. 
    \end{example}
    
    The order of a graph $|G|$ is the number of vertices of $G$, while the size of a graph $\operatorname{size}(G)$ is the number of edges of $G$.
    For example, for a tree $T$ we have $|T|=\operatorname{size}(T)+1$.
    If one labels the vertices of a sandpile $G$ as $v_1,\dots,v_{|G|}$, where the graph $G$ may be weighted with positive integer weights, then the \textbf{Laplacian} of $G$ is the $|G| \times |G|$ matrix $\Delta = D^\mathrm{T} - A^\mathrm{T}$, where $D$ is the diagonal matrix where $D_{ii} = d_{v_i}$, and $A$ is the adjacency matrix of $G$.
    That is, if~$a_{ij}$ is the weight of the edge from vertex $v_i$ to $v_j$, and $d_i$ is the degree of~$v_i$,
    \begin{equation*}
        \Delta_{ij} = \begin{cases}
        -a_{ij} & \text{for } i \neq j, \\
        d_i & \text{for } i = j.
        \end{cases}
    \end{equation*}
    
    
    The \textbf{reduced Laplacian} $\Delta '$ of $G$ is obtained by removing from $\Delta$ the row and column corresponding to the sink.
    To explicitly refer to the reduced Laplacian with the sink at $v_i$, the notation~$\Delta^{(i)}$ is used.
    We can represent the firing of a non-sink vertex $v$ as the subtraction of the column of $\Delta'$ corresponding to $v$ from the chip configuration. 
    
    \begin{example}\label{example: Laplacian C4}
    The cycle graph on four vertices $C_4$ has Laplacian
    \begin{equation*}
        \Delta = \begin{bmatrix}
        2  & -1 & 0  & -1 \\
        -1 & 2  & -1 & 0  \\
        0  & -1 & 2  & -1 \\
        -1 & 0  & -1 & 2  \\
        \end{bmatrix}
    \end{equation*}
    and the reduced Laplacian, for sink $v_1$, is
    \begin{equation*}
        \Delta' = \begin{bmatrix}
        2 & -1 & 0 \\
        -1 & 2 & -1 \\
        0 & -1 & 2 \\
        \end{bmatrix}.
    \end{equation*}
    \end{example}

    In order to view configurations before and after firing as equivalent, we define the \textbf{sandpile group} of $G$ to be the group quotient
    
    \begin{equation*}
        \mathcal{S}(G)=\mathbb{Z}^{|G|-1} / \Delta^{\prime} \mathbb{Z}^{|G|-1}.
    \end{equation*}
    
    From the definition of the sandpile group, we see that $\left|\mathcal{S}(G)\right|=\left|\Delta'\right|$, where $\left|\mathcal{S}(G)\right|$ denotes the order of the sandpile group $\mathcal{S}(G)$ and $\left|\Delta'\right|$ denotes the determinant of the reduced Laplacian.
    This holds regardless of the choice of sink.
    By the Matrix-Tree theorem (see for example \cite{stanley2013algebraic}), $\left|\Delta'\right|$ is the number of spanning trees of $G$.
    It is also shown in \cite{biggs1999chip} that each equivalence class of $\mathcal{S}(G)$ contains exactly one recurrent configuration, and that the recurrent configurations form an abelian group with the operation being defined as the stabilization of the sum of two recurrent configurations.
    
    
    A configuration $c$ over graph $G$ with $n$ vertices is \textbf{equivalent} to another configuration $d$ when they have the same image in $\mathcal{S}(G)$, meaning they lie in the same equivalence class in $\mathcal{S}(G)$, or in other words,
    \begin{equation*}
        c\equiv d \Leftrightarrow \text{there exists } \, \mathbf{v} \in \mathbb{Z}^{|G|-1}: c - d = \Delta' \mathbf{v}.
    \end{equation*}

    Notice that as $\Delta'$ is non-singular by the Matrix-Tree theorem, the vector $(\Delta')^{-1} (c-d)$ is the unique solution to the equation $c = d + \Delta' \mathbf{v}$, and thus
    \begin{equation*}
        c \equiv d \Longleftrightarrow (\Delta')^{-1} (c-d) \in \mathbb{Z}^{|G|-1}.
    \end{equation*}
    
    However, two configurations being equivalent does not necessarily imply there exists a firing sequence that takes one to the other.
    Hence, we introduce the ideas of backfiring and unrestricted firing. 
    To fire a vertex $v$ is equivalent to subtracting a column of $\Delta'$ corresponding to $v$, and similarly to \textbf{backfire} $v$ is equivalent to \textit{adding} that column.
    Backfiring has been studied in \cite{dhar1994inverse} and \cite{caracciolo2012multiple} under the names ``untoppling" and ``anti-toppling", respectively.
    \textbf{Unrestricted firing} is when vertices are allowed to fire and/or backfire regardless of the number of chips they have.
    Furthermore, the sink is allowed to fire, where firing the sink corresponds to backfiring all non-sink vertices.
    
    Of interest in the sandpile group is the identity element, which leads us to the definition of the recurrent identity.
    
    \begin{definition}\label{definition: recurrent identity}
        The \textbf{recurrent identity} is the recurrent configuration equivalent to the all-zero configuration.
    \end{definition}
    
    Recall that each equivalence class contains exactly one recurrent configuration, and thus the recurrent identity is well-defined.
    

    \begin{example}
    Let $G=K_3$. We look at the sandpile of $G$ with arbitrary sink $s$. We will show that the recurrent identity $c$ is equivalent to the all-zero configuration $d$, where our two chip configurations are
    \begin{equation*}
        c = \begin{bmatrix} 1 \\ 1 \end{bmatrix}, \quad d = \begin{bmatrix} 0 \\ 0 \end{bmatrix}
    \end{equation*}
    and the reduced Laplacian is
    \begin{equation*}
        \Delta' = \begin{bmatrix} 2 & -1 \\ -1 & 2 \end{bmatrix}.
    \end{equation*}
    Now we notice that
    \begin{equation*}
        (\Delta')^{-1} (c-d) = \begin{bmatrix}[1.3] \frac{2}{3} & \frac{1}{3} \\ \frac{1}{3} & \frac{2}{3} \end{bmatrix} \begin{bmatrix} 1 \\ 1 \end{bmatrix} = \begin{bmatrix} 1 \\ 1 \end{bmatrix} \in \mathbb{Z}^2.
    \end{equation*}
    
    Hence, $c\equiv d$, and $d$ can be obtained from $c$ by unrestrictedly firing (see \cref{section:CMIP} for a formal definition) the two non-sink vertices once each. See \cref{fig:equivalence} for a visual example of the firing sequence taking $c$ to $d$, and thus the equivalence of $c$ and $d$.
    
    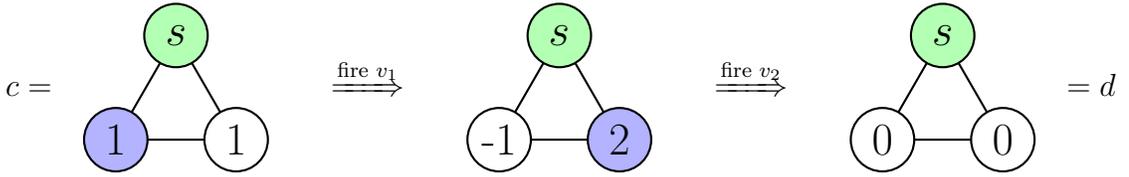
\begin{figure}[htbp!]
        \centering
        $c=$
        \begin{subfigure}{0.2\textwidth}
            \centering
            \begin{tikzpicture}[scale=0.8,font=\Large,baseline]
                \draw[thick] (-1,-0.866) -- (1,-0.866) -- (0,0.866) -- (-1,-0.866);
                \filldraw[color=black,fill=green!30,thick] (0,0.866) circle (15pt);
                \filldraw[color=black,fill=blue!30,thick] (-1,-0.866) circle (15pt);
                \filldraw[color=black,fill=white,thick] (1,-0.866) circle (15pt);
                \node at (0,0.866) {$s$};
                \node at (-1,-0.866) {1};
                \node at (1,-0.866) {1};
            \end{tikzpicture}
        \end{subfigure}
        \quad
        $\xRightarrow{\text{fire } v_1}$
        \quad
        \begin{subfigure}{0.2\textwidth}
            \centering
            \begin{tikzpicture}[scale=0.8,font=\Large,baseline]
                \draw[thick] (-1,-0.866) -- (1,-0.866) -- (0,0.866) -- (-1,-0.866);
                \filldraw[color=black,fill=green!30,thick] (0,0.866) circle (15pt);
                \filldraw[color=black,fill=white,thick] (-1,-0.866) circle (15pt);
                \filldraw[color=black,fill=blue!30,thick] (1,-0.866) circle (15pt);
                \node at (0,0.866) {$s$};
                \node at (-1,-0.866) {-1};
                \node at (1,-0.866) {2};
            \end{tikzpicture}
        \end{subfigure}
        \quad
        $\xRightarrow{\text{fire } v_2}$
        \quad
        \begin{subfigure}{0.2\textwidth}
            \centering
            \begin{tikzpicture}[scale=0.8,font=\Large,baseline]
                \draw[thick] (-1,-0.866) -- (1,-0.866) -- (0,0.866) -- (-1,-0.866);
                \filldraw[color=black,fill=green!30,thick] (0,0.866) circle (15pt);
                \filldraw[color=black,fill=white,thick] (-1,-0.866) circle (15pt);
                \filldraw[color=black,fill=white,thick] (1,-0.866) circle (15pt);
                \node at (0,0.866) {$s$};
                \node at (-1,-0.866) {0};
                \node at (1,-0.866) {0};
            \end{tikzpicture}
        \end{subfigure}
        $=d$
        
        \caption{Equivalence of the recurrent identity and the all-zero configuration. The sink is colored green, and the vertex being fired is colored blue.}\label{fig:equivalence}

    \end{figure}
    \end{example}
    
\begin{figure}[htbp!]
        \centering
        \begin{subfigure}{0.4\textwidth}
            \centering
            \includegraphics[width=\textwidth]{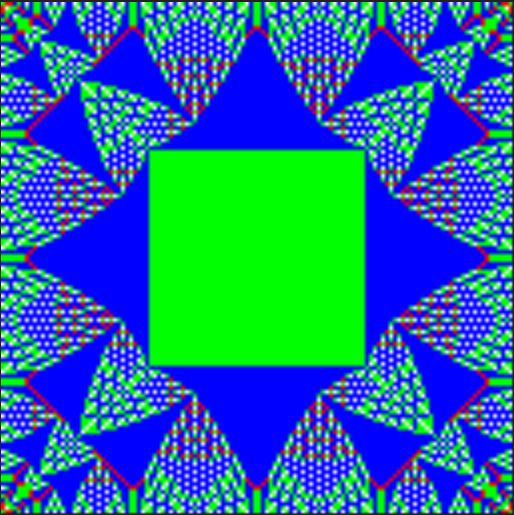}
            \caption{$128 \times 128$ square grid}
            \label{subfigure: Identity 128}
        \end{subfigure}
        \quad
        \begin{subfigure}{0.4\textwidth}
            \centering
            \includegraphics[width=\textwidth]{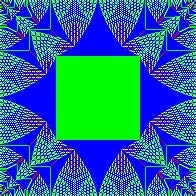}
            \caption{$198 \times 198$ square grid}
            \label{subfigure: Identity 198}
        \end{subfigure}
        \caption{The recurrent identity of the $128\times128$ and $198\times198$ square grids, as shown in \cite{holroyd2008chip}. 0,~1,~2, and 3 chips are displayed as orange, red, green, and blue, respectively.}
        \label{figure: Identity Square Grid}
    \end{figure}

    The recurrent identity has been studied previously, including some limiting behavior for $\Z^2$ lattice graphs as in \cite{caracciolo2008explicit} and \cite{le2002identity}.
    For example, the recurrent identity for the $128 \times 128$ and $198 \times 198$ square grids with a boundary sink, or a sink connected to all the boundary vertices, is shown in \cref{figure: Identity Square Grid}.
    Many intriguing questions about the identity are still open.
    
\section{The Complete Maximal Identity Property}\label{section:CMIP}

    The maximal stable configuration is guaranteed to be recurrent over all sandpiles, even when there is only one recurrent configuration, and is by far the simplest recurrent configuration.
    Furthermore, of interest among the recurrent elements is the identity of the abelian group, the recurrent identity.
    Hence, it is of interest to know for which graphs the maximal stable configuration is the recurrent identity, leading to our newly created definition of the maximal identity property.
    
    \begin{definition}\label{definition:MIP}
        A graph $G$ is said to have the \textbf{maximal identity property} at a vertex $s$ if, having chosen $s$ as the sink, the maximal stable configuration is the recurrent identity of the sandpile $G$. 
    \end{definition}
    

    
    The following lemma provides an equivalent condition for the maximal identity property.
    
    
    \begin{lemma}\label{lemma: unrestricted firing MIP}
         The sandpile $G$ has the maximal identity property if and only if there exists an unrestricted firing sequence that takes $m_G$ to $\mathbf{0}$.
    \end{lemma}
    \begin{proof}
        The maximal stable configuration $m_G$ is always recurrent.
        An unrestricted firing sequence between $m_G$ and $\mathbf{0}$ exists if and only if the two configurations are equivalent, which occurs if and only if $m_G$ is the recurrent identity, resulting in the maximal identity property.
    \end{proof}
    
    Another necessary and sufficient condition for a sandpile $G$ to have the maximal identity property is for $\operatorname{Stab}(m_G+m_G)=m_G$.
    
    The complete maximal identity property, developed for the first time in this paper, allows any vertex of a graph to be chosen as the sink while preserving the maximal identity property.
    
    \begin{definition}\label{definition:CMIP}
    A graph is said to have the \textbf{complete maximal identity property} if for all vertices $v$, it has the maximal identity property for sink $v$.
    \end{definition}

\begin{figure}[h!]
        \centering
        \begin{subfigure}{0.3\textwidth}
            \centering
            \begin{tikzpicture}[thick,scale=0.9]
                \draw (0,2) -- (0.951*2,0.309*2) -- (0.588*2,-0.809*2) -- (-0.588*2,-0.809*2) -- (-0.951*2,0.309*2) -- (0,2);
                \draw (0,1) -- (0,2);
                \draw (0.951,0.309) -- (0.951*2,0.309*2);
                \draw (0.588,-0.809) -- (0.588*2,-0.809*2);
                \draw (-0.588,-0.809) -- (-0.588*2,-0.809*2);
                \draw (-0.951,0.309) -- (-0.951*2,0.309*2);
                \draw (0,1) -- (-0.588,-0.809) -- (0.951,0.309) -- (-0.951,0.309) -- (0.588,-0.809) -- (0,1);
                \foreach \x in {1,2} {
                    \filldraw[color=black,fill=white] (0 * \x, 1 * \x) circle (7pt);
                    \foreach \y in {1,-1} {
                        \filldraw[color=black,fill=white] (0.951 * \x * \y, 0.309 * \x) circle (7pt);
                        \filldraw[color=black,fill=white] (0.588 * \x * \y, -0.809 * \x) circle (7pt);
                    }
                }
            \end{tikzpicture}
            \caption{The Petersen graph}
            \label{subfigure: Petersen graph}
        \end{subfigure}
        \quad
        \begin{subfigure}{0.3\textwidth}
            \centering
            \begin{tikzpicture}[thick,scale=1.2]
                \draw (-0.5,-0.866) -- (1.5,-0.866) -- (2,0) -- (1,1.732) -- (0,1.732) -- (-1,0) -- (-0.5,-0.866) -- (0,0) -- (-0.5,0.866);
                \draw (0,0) -- (-1,0);
                \draw (0,1.732) -- (0.5,0.866) -- (1.5,0.866);
                \draw (0.5,0.866) -- (1,1.732);
                \draw (0.5,-0.866) -- (1,0) -- (2,0);
                \draw (1,0) -- (1.5,-0.866);
                \filldraw[color=black,fill=white] (0,0) circle (6pt);
                \filldraw[color=black,fill=white] (-0.5,.866) circle (6pt);
                \filldraw[color=black,fill=white] (1,1.732) circle (6pt);
                \filldraw[color=black,fill=white] (0,1.732) circle (6pt);
                \filldraw[color=black,fill=white] (1.5,-.866) circle (6pt);
                \filldraw[color=black,fill=white] (2,0) circle (6pt);
                \filldraw[color=black,fill=white] (1.5,0.866) circle (6pt);
                \filldraw[color=black,fill=white] (0.5,0.866) circle (6pt);
                \filldraw[color=black,fill=white] (0.5,-0.866) circle (6pt);
                \filldraw[color=black,fill=white] (1,0) circle (6pt);
                \filldraw[color=black,fill=white] (-0.5,-0.866) circle (6pt);
                \filldraw[color=black,fill=white] (-1,0) circle (6pt);
            \end{tikzpicture}
            \caption{The 3-diamond ring}
            \label{subfigure: 3-diamond Ring}
        \end{subfigure}
        \quad
        \begin{subfigure}{0.3\textwidth}
            \centering
            \begin{tikzpicture}[thick,scale=1.2]
                \draw (1.5,0) -- (0,1.5) -- (-1.5,0) -- (0,-1.5) -- (1.5,0);
                \draw (-1.5,0) -- (-0.5,0);
                \draw (0.5,0) -- (1.5,0);
                \draw (0,-1.5) -- (0,-0.5);
                \draw (0,0.5) -- (0,1.5);
                \draw (1,-0.5) -- (0.5,0) -- (1,0.5);
                \draw (-0.5,1) -- (0,0.5) -- (0.5,1);
                \draw (-1,0.5) -- (-0.5,0) -- (-1,-0.5);
                \draw (-0.5,-1) -- (0,-0.5) -- (0.5,-1);
                \filldraw[color=black,fill=white] (0,-1.5) circle (5pt);
                \filldraw[color=black,fill=white] (0.5,1) circle (5pt);
                \filldraw[color=black,fill=white] (0.5,0) circle (5pt);
                \filldraw[color=black,fill=white] (-1,0.5) circle (5pt);
                \filldraw[color=black,fill=white] (0,1.5) circle (5pt);
                \filldraw[color=black,fill=white] (1.5,0) circle (5pt);
                \filldraw[color=black,fill=white] (-0.5,0) circle (5pt);
                \filldraw[color=black,fill=white] (1,-0.5) circle (5pt);
                \filldraw[color=black,fill=white] (-1.5,0) circle (5pt);
                \filldraw[color=black,fill=white] (0.5,-1) circle (5pt);
                \filldraw[color=black,fill=white] (-1,-0.5) circle (5pt);
                \filldraw[color=black,fill=white] (-0.5,-1) circle (5pt);
                \filldraw[color=black,fill=white] (0,0.5) circle (5pt);
                \filldraw[color=black,fill=white] (-0.5,1) circle (5pt);
                \filldraw[color=black,fill=white] (0,-0.5) circle (5pt);
                \filldraw[color=black,fill=white] (1,0.5) circle (5pt);
            \end{tikzpicture}
            \caption{The 4-diamond ring}
            \label{subfigure: 4-diamond Ring}
        \end{subfigure}
        \caption{The Petersen graph and the $n$-diamond ring for $n=3,4$}
        \label{figure: Petersen and Diamong Ring}
    \end{figure}
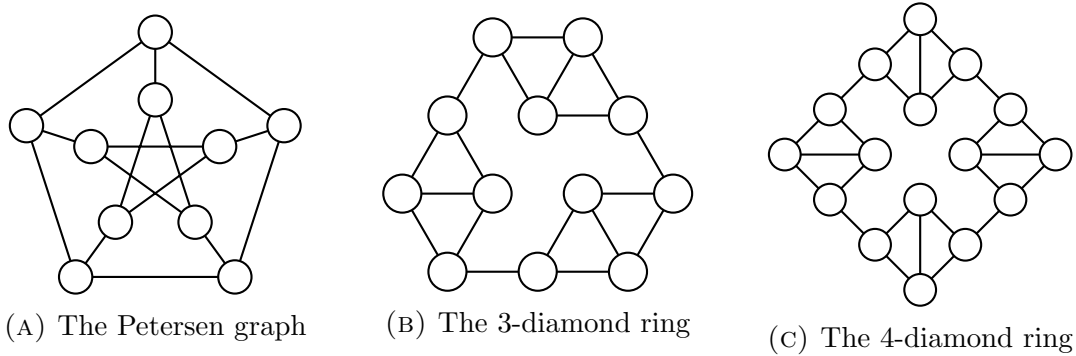

    \begin{example}\label{example: CMIP Petersen and Diamond Ring}
        Both the Petersen graph and the $n$-diamond ring for all positive integers $n$ have the complete maximal identity property (see \cref{figure: Petersen and Diamong Ring} for the Petersen graph and the $n$-diamond ring for $n=3,4$). The proofs are straightforward but tedious, and thus not included in this paper.
    \end{example}
    
    \begin{proposition}\label{proposition: CMIP complete odd cycle tree}
        All complete graphs $K_n$, odd cycles $C_{2n+1}$, and trees have the complete maximal identity property.
    \end{proposition}
    \begin{proof}
        \textbf{Case 1:} $K_n$.
        By symmetry, we only need to prove that $m_{K_n} \equiv \mathbf{0}$ for one sink of $K_n$.
        We prove this via the unrestricted firing method of \cref{lemma: unrestricted firing MIP}.
        The configuration~$m_{K_n}$ has~$n-2$ chips at each non-sink vertex. Backfire the sink $n-2$ times to yield the all-zero configuration.
        
        
        \textbf{Case 2:} $C_{2n+1}$.
        By symmetry, we only need to prove that $m_{C_{2n+1}} \equiv \mathbf{0}$ for one sink of $C_{2n+1}$.
        We prove this via the unrestricted firing method.
        Note that $m_{C_{2n+1}} \equiv \mathbf{1}$.
        Define the \textbf{rank} of a vertex to be the edge length of the shortest path from the sink to that vertex.
        Fire both vertices of rank $n$ once. This results in a net movement on one chip to each of the vertices of rank $n-1$, ridding the vertices of rank $n$ of chips.
        Now, as each vertex of rank $n-1$ has 2 chips, fire all vertices with rank at least $n-1$ twice to eliminate all chips from those vertices.
        Continue in this manner, firing all vertices with rank $n-k$ or higher $k+1$ times.
        This results in the all-zero configuration.
        
        \textbf{Case 3:} Trees.
        There is only one spanning tree of a tree, the tree itself.
        Hence, we apply the Matrix Tree Theorem to find that there is only one equivalence class in the sandpile group, regardless of the choice of sink.
        Thus, there is only one recurrent element; as $m_G$ is always recurrent, it is the only recurrent configuration.
        Therefore, $m_G$ must be the recurrent identity.
    \end{proof}
    
    The fact that a tree has a trivial sandpile group gives it a unique effect when added to other graphs.
    The concept of adding a tree of size $n$ to a vertex $v\in G$ is simply taking any tree~$T$ on~$n+1$ vertices and $n$ edges, then taking the disjoint union of $T$ and $G$ by combining their vertex sets together and their edge sets together, and finally merging a vertex of~$T$ with vertex $v$.
    We will show in \cref{lemma:addtreeCMIPidentityconfiguration} that ultimately only the size of the tree matters.
    See \cref{Diamond} and \cref{Diamond+} for an example on how trees are added to graphs.
    Notice that adding the tree of size 0, a single vertex, does not change the graph.
    The following theorem states how one can attach trees to any connected graph in the aforementioned manner to result in a graph that has the complete maximal identity property.
    
    \begin{figure}[htbp!]
    \centering
    \begin{minipage}{0.45\textwidth}
        \centering
        \begin{tikzpicture}[scale=0.8,font=\Large,baseline]
            \draw[thick] (-1,-1) -- (1,-1) -- (1,1) -- (-1,1) -- (-1,-1) -- (1,1);
            \filldraw[color=black,fill=white,thick] (-1,-1) circle (15pt);
            \filldraw[color=black,fill=white,thick] (-1,1) circle (15pt);
            \filldraw[color=black,fill=white,thick] (1,1) circle (15pt);
            \filldraw[color=black,fill=white,thick] (1,-1) circle (15pt);
            \node at (-1,-1) {$v_1$};
            \node at (-1,1) {$v_2$};
            \node at (1,1) {$v_3$};
            \node at (1,-1) {$v_4$};
        \end{tikzpicture}
        \caption{The Diamond Graph}\label{Diamond}
    \end{minipage}
    \quad
    \begin{minipage}{0.5\textwidth}
        \centering
        \begin{tikzpicture}[scale=0.8,font=\Large,baseline]
            \draw[thick] (-1,-1) -- (1,-1) -- (1,1) -- (-1,1) -- (-1,-1) -- (1,1);
            \draw[thick] (1,-1) -- (2.414,-2.414);
            \filldraw[color=black,fill=white,thick] (-1,-1) circle (15pt);
            \filldraw[color=black,fill=white,thick] (-1,1) circle (15pt);
            \filldraw[color=black,fill=white,thick] (1,1) circle (15pt);
            \filldraw[color=black,fill=white,thick] (1,-1) circle (15pt);
            \filldraw[color=black,fill=white,thick] (2.414,-2.414) circle (15pt);
            \node at (-1,-1) {$v_1$};
            \node at (-1,1) {$v_2$};
            \node at (1,1) {$v_3$};
            \node at (1,-1) {$v_4$};
            \node at (2.414,-2.414) {$v_5$};
        \end{tikzpicture}
        \caption{The Diamond Graph with a tree of size 1 added to vertex $v_4$. This graph has the maximal identity property for sink at $v_2$.}\label{Diamond+}
    \end{minipage}
    \end{figure}


    \begin{theorem}\label{theorem: Add tree CMIP}
        Given any connected graph $G$, there exists infinitely many graphs derived from adding trees to $G$ that have the complete maximal identity property. 
    \end{theorem}

    To prove \cref{theorem: Add tree CMIP}, we need the following lemma.

    \begin{lemma}\label{lemma:addtreeCMIPidentityconfiguration}
        Let $G$ be any graph. Let $\{T_v\}$ be a family of trees labeled by vertices of $G$ and let~$G'$ be the graph obtained from $G$ by adding each tree $T_v$ to vertex $v$.
        Then $G'$ has the complete maximal identity property if and only if the configuration $c$ over $G$ with $d_v + |T_v| - 2$ chips at all vertices $v$ is equivalent to the identity for all selections of sinks, where $|T_v|$ is the number of vertices in the tree.
    \end{lemma}
    \begin{proof}
        We first prove the if direction.
        Pick a sink $s\in G \subseteq G'$.
        Let $c$ be the chip configuration in the lemma statement.
        We observe that each tree $T_v$ has size $c_v - (m_G)_v$.
        As trees have a trivial sandpile group, there exists an unrestricted firing sequence that moves all the chips of $m_{G'}$ in $T_v$ to $v$.
        Each edge of $T_v$ provides one additional chip to $m_{G'}$ from $m_G$, and thus once all the chips are moved towards their attachment point $v$ in $G$, we find that the configuration that results is configuration~$c\equiv m_{G'}$.
        As~$c\equiv \mathbf{0_G}$, there exists an unrestricted firing sequence that takes $c$ to $\mathbf{0}$ in~$G$.
        Use this same sequence on the chip configuration $c$ over $G'$, except whenever vertex $v$ in $G$ is fired, fire vertex~$v$ and all the vertices of $T_v$ in $G'$.
        This results in $\mathbf{0_{G'}}$, thus illustrating that $G'$ has the maximal identity property for sink $s$.
        The choice of $s$ was arbitrary, so it applies for all vertices in $G$.
        
        
        However, $G'$ also has the vertices in the trees that were added.
        We must prove that $G'$ has the maximal identity property for these vertices as sinks.
        Pick a sink $s\in G'$ from tree $T_v$ for arbitrary vertex $v\in G$.
        Use the unrestricted firing sequence that took $m_{G'}$ to $\mathbf{0}$ for sink $v$, to result in a configuration equivalent to $m_G$ that has all of the chips at vertex $v$.
        As all trees have a trivial sandpile group, there exists a unrestricted firing sequence that takes all the chips at vertex $v$ to $s$, if only $T_v$ is considered. Now, whenever vertex $v$ needs to be fired, fire all vertices in $G$, and all vertices of $T_{v'}$ for all $v'\neq v$. This allows all firings in $T_v\subset G'$ to work inside $G'$. This demonstrates that $m_{G'}\equiv \mathbf{0}$, and thus that $G'$ has the complete maximal identity property.
        
        To prove the only if direction, assume $G'$ has the complete maximal identity property. Reversing the unrestricted firing sequences used to move all the chips of $m_{G'}$ in $T_v$ to $v$, we find $c \equiv m_{G'}$ as the $|T_v|-1$ chips $T_v$ needs for the maximal stable configuration are taken away from $v$ to result in $v$ having $d_v - 1$ chips, the number of chips it needs for the maximal stable configuration. As $G'$ has the complete maximal identity property, $c \equiv m_{G'} \equiv \mathbf{0_{G'}}$. 
    \end{proof}
    
    We now can prove \cref{theorem: Add tree CMIP}.
    
    \begin{proof}[Proof of \cref{theorem: Add tree CMIP}]
        Say we have a graph $G$.
        By the Matrix-Tree Theorem, regardless of the choice of sink, $|\Delta'|$ is a constant integer.
        Say $|\Delta'|=k$. By \cref{lemma:addtreeCMIPidentityconfiguration}, to prove the result it suffices to find infinitely many configurations $c$, where $c$ is over all vertices, such that $c \geq m_G$ and for all selections of sinks, $c\equiv\mathbf{0}$.
        We will create such a configuration as follows.
        For each vertex $v$, let the number of chips on it be equal to a multiple of $k$ that is greater than or equal to $d_v-1$.
        This configuration can be represented by a vector $k\mathbf{x}$, where $\mathbf{x}$ is a vector composed of nonnegative integers.
        For a particular selection of sink, to prove that this configuration is equivalent to the identity, we must show that there exists a vector $\mathbf{y}$ with integer entries such that
        \begin{equation*}
            \Delta' \mathbf{y} = k\mathbf{x}.
        \end{equation*}
        But we know that $|\Delta'|=k$, so using the fact that the adjoint matrix of an integer matrix (which~$\Delta'$ is) has integer entries, we know that $k(\Delta')^{-1}$ has integer entries. Hence
        \begin{equation*}
            \mathbf{y} = k(\Delta')^{-1} \mathbf{x},
        \end{equation*}
        and as $\mathbf{x}$ has integer entries, $\mathbf{y}$ has integer entries.
        As there are infinitely many multiples of a positive integer greater than a fixed value, there are infinitely many valid configurations $c$.
        Thus, there exists infinitely many graphs consisting of the given graph with trees attached to it that have the complete maximal identity property.
    \end{proof}
    
    Because we may add trees to graphs to give them the complete maximal identity property, we wish to have a notion of irreducibility that eliminates such graphs which have trees added to them.
    This leads us to the classical notion of a biconnected graph.
    
    \begin{definition}\label{definition: biconnected}
        A graph is \textbf{biconnected} if it remains connected even if one removes any single vertex and its incident edges.
    \end{definition}
    
    In other words, for any two vertices in a biconnected graph, there exist at least two vertex-disjoint paths that connect them.
    
    A search over all biconnected graphs with 11 or fewer vertices found that asides from odd cycles and complete graphs, there were only three and two biconnected graphs (up to isomorphism) with 8 and 10 vertices, respectively.
    These include the 2-diamond ring and the Petersen graph from \cref{example: CMIP Petersen and Diamond Ring}.
    In addition, there are three other biconnected graphs with 12 vertices known to have the complete maximal identity property, including the 3-diamond ring; however, the search over all biconnected graphs with 12 vertices is too computationally intensive to complete.
    
    The observed lack of any other biconnected graphs with an odd number of vertices that possess the complete maximal identity property prompts the following question.
    
    \begin{question}\label{question: odd biconnected CMIP}
        Are the only biconnected graphs with an odd number of vertices that possess the complete maximal identity property cycle graphs and complete graphs?
    \end{question}

\section{Necessary Conditions for the Complete Maximal Identity Property}\label{section: necessary and sufficient}



In order to create necessary conditions for the complete maximal identity property, we first make the following definition.

\begin{definition}\label{definition:compatible}
A vector $c\in\Z^n$ is \textbf{compatible} if for any sink $v_i$ we have
\begin{equation*}
    (c_1,\ldots,c_{i-1},c_{i+1},\ldots,c_n)\in\Delta^{(i)}\Z^{n-1},
\end{equation*}
where the vector on the left hand side will be denoted by $c^{(i)}$.
An equivalent definition is that for any choice of sink $v_i$,
\begin{equation*}
    (\Delta^{(i)})^{-1} c^{(i)} \in \Z^{|G|-1}.
\end{equation*}

\end{definition}

Note that if $a$ and $b$ are compatible, then $a\pm b$ is compatible.

The following lemma considers compatible configurations where only a single vertex contains chips.
We use $\mathbf{e}_i$ to denote the $i$th standard basis vector with a 1 in the $i$th coordinate and 0s elsewhere.
\begin{lemma}\label{lemma:compatible sum entries}
If $c$ is compatible and $s=c_1 + \cdots + c_n$, then $s\mathbf{e}_i$ is compatible for all $i$.
\end{lemma}
\begin{proof}
Notice that if $c\in\Z^n$ is in the image of $\Delta$, then $c$ is compatible.
Since $(c_1,\dots,c_n)$ is compatible, then
\[    (c_1,\dots,c_{i-1}, c_{i+1}, \dots, c_n)=\sum_{j\in\{1,\dots,n\}\setminus\{i\}}\alpha_j \Delta^{(i)}_j, \]
where $\alpha_j\in\Z$ and $\Delta^{(i)}_j$ denotes the column of $\Delta^{(i)}$ corresponding to vertex $v_j$. Then 
\[    \left(c_1,\dots,c_{i-1},-\sum_{j\in\{1,\dots,n\}\setminus\{i\}} c_j,c_{i+1},\dots,c_n\right)=\sum_{j\in\{1,\dots,n\}\setminus\{i\}}\alpha_j \Delta_j    \]
is in the image of $\Delta$ and thus is compatible.
Subtracting the two compatible vectors, then $s\mathbf{e}_i$ is compatible for any $i$ where $s = c_1 + \cdots +c_n$.
\end{proof}


Clearly, $I=\{d\mid d\mathbf{e}_i\text{ is compatible for any }i\}$
forms an ideal in $\Z$.
For a graph $G$ with the complete maximal identity property, we know $k\in I$, where $k=|\Delta'|$, and
\[ \sum_{v\in G} (\mathrm{deg}(v)-1) = 2 \cdot \operatorname{size}(G) - |G| \in I, \]
where $\operatorname{size}(G)$ denotes the number of edges of $G$, so their greatest common divisor is in $I$.

As $\Z$ is a principal ideal domain, for each graph $G$, $I=(x)$ for some positive integer $x$. Notice that $x$ is the smallest positive element of $I$. This leads us to our definition of the minimal compatibility number of a graph.

\begin{definition}\label{definition: minimal compatibility number}
    The \textbf{minimal compatibility number} of a graph is the positive integer $x$ such that $(x) = I := \{d\mid d\mathbf{e}_i\text{ is compatible for any }i\}$.
\end{definition}

\begin{lemma}\label{lemma: minimal compatibility number = 1 >>> tree}
    The minimal compatibility number of a graph $G$ is $1$ if and only if $G$ is a tree.
\end{lemma}
\begin{proof}
We first prove the only if direction.
If the minimal compatibility number of a graph is 1, then every chip configuration is in the integer image of the reduced Laplacian and thus equivalent to each other.
Hence, the sandpile group $\mathcal{S}(G)$ is the trivial group, so by the Matrix-Tree Theorem $G$ must have only 1 spanning tree, or in other words, is a tree itself.

For the if direction, a tree has a trivial sandpile group and thus every chip configuration is equivalent to each other, meaning that each chip configuration is in the integer image of the reduced Laplacian, and thus the minimal compatibility number of the tree is 1.
\end{proof}

Notice that $d\mathbf{e}_i$ is compatible if each least common denominator of the column of $(\Delta')^{-1}$ corresponding to vertex $i$ for all sinks $s$ divides $d$.
Hence, $d\in I$ if the least common denominator of all entries of the inverses of all reduced Laplacians of a graph $G$ divides $d$, and thus the minimal compatibility number of a graph is the least common denominator of all entries of the inverses of all reduced Laplacians of a graph $G$.

\begin{proposition}\label{proposition: CMIP >>> gcd k sum(msc) > 1}
    If a non-tree graph $G$ has the complete maximal identity property, then
    \[\gcd\left(\left|\Delta'\right|,2\cdot\operatorname{size}(G)-|G|\right)>1. \]
\end{proposition}
\begin{proof}
    We will prove the contrapositive.
    
    If for a non-tree graph $G$, 
    \[ \gcd\left(\left|\Delta'\right|,2\cdot\operatorname{size}(G)-|G|\right)=1, \]
    the graph does not have the complete maximal identity property, as if it did then 
    \[ \gcd\left(\left|\Delta'\right|,2\cdot\operatorname{size}(G)-|G|\right)=1\in I, \]
    so the minimal compatibility number of $G$ is 1, which by \cref{lemma: minimal compatibility number = 1 >>> tree} means $G$ is a tree, yielding a contradiction.
\end{proof}

\begin{proposition}\label{proposition: CMIP >>> x < n^2}
    If a graph $G$ has the complete maximal identity property, then the minimal compatibility number $x$ satisfies $x \leq |G|^2 - 2|G|$ for $|G| > 2$, and $x = 1$ if $|G| = 2$.
\end{proposition}
\begin{proof}
    If $|G| = 2$, then $G$ must be the connected graph with two vertices, which is a tree and thus the result follows from \cref{lemma: minimal compatibility number = 1 >>> tree}.
    
    We now assume $|G| > 2$.
    As
    \[ 2\cdot \operatorname{size}(G) - |G| \in I, \]
    we have
    \[ x \leq 2\cdot\operatorname{size}(G) - |G|, \]
    as $x$ is the smallest positive element of $I$ and $2\cdot \operatorname{size}(G) - |G| > 0$.
    This is because $G$ is a connected graph so $\operatorname{size}(G) \geq |G| - 1$, and thus $2\cdot \operatorname{size}(G) - |G| \geq |G| - 2 > 0$.
    
    The maximum value of $\operatorname{size}(G)$ of a graph with fixed order is $\frac{|G|(|G|-1)}{2}$, so
    \[ x \leq 2\cdot \operatorname{size}(G) - |G| \leq \left(|G|^2 - |G|\right) - |G| = |G|^2 - 2|G|. \]
\end{proof}

\begin{proposition}\label{proposition: x Kn Cn = n}
    The minimal compatibility number for both $K_n$ and $C_n$ is $n$ for $n \geq 3$.
\end{proposition}
\begin{proof}
    We will show $n\mathbf{e}_i$ is compatible, with an irreducible firing vector (i.e. the elements have no common factor).
    The \textbf{firing vector} for a configuration $c$ is the vector $(\Delta')^{-1} c$, where each element of the vector corresponds to how many times the corresponding non-sink vertex must fire from the all-zero configuration to reach $c$.
    Note that the firing vector may not necessarily have integer entries.
    As the reduced Laplacian is non-singular by the Matrix-Tree Theorem, the firing vector is well-defined, existing uniquely as the firing vector that takes the all-zero configuration to any configuration $c$.
    If $n\mathbf{e}_i$ is compatible, but the firing vector is irreducible, meaning all the entries are integers which do not share a nontrivial common factor, then~$n'\mathbf{e}_i$ for all positive integers $n' < n$ will not have a firing vector that has all integer entries;
    this is because that firing vector can be obtained by multiplying the firing vector for $n\mathbf{e}_i$ by $\frac{n'}{n}$.
    
    For the complete graph $K_n$, fire the sink and then backfire $v_i$ to result in $n$ chips at $v_i$.
    The firing vector consists of $-1$ at all non-sink vertices except $v_i$ which is $-2$.
    Hence, this is an irreducible firing vector, and $x=n$.
    
    For the cycle graph $C_n$, we will first show that $n\mathbf{e}_i$ is compatible for all $i$.
    Notice that on the cycle graph, one may select a proper connected subgraph (i.e. an ``arc" of the circle) and fire all of them, resulting in the outer vertices losing one vertex and the vertices adjacent to them but not in the arc gaining a chip.
    This may be repeated for iteratively larger arcs, each time including one more vertex on each end.
    Doing so enables us to send those chips an arbitrary distance away from the ends of the original arc (as long as there are no self-intersection issues).
    
    We will first look at odd $n$.
    From the sink, give each vertex a single chip by pairing vertices equidistant from the sink, and having the sink fire chips to those two vertices.
    Then, pair vertices equidistant from vertex $v_i$ and fire their chips to vertex $v_i$.
    This results in vertex $v_i$ having all $n$ chips, with no other vertex having chips.
    See \cref{figure: minimal compatibility example C7} for an example of the described firing process.
    
    \begin{figure}[h!]
        \centering
        \includegraphics[width=0.2\textwidth,valign=c]{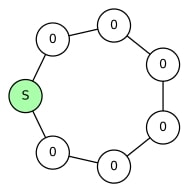}
        $\Rightarrow$
        \includegraphics[width=0.2\textwidth,valign=c]{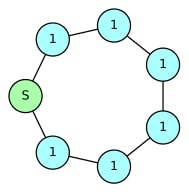}
        $\Rightarrow$
        \includegraphics[width=0.2\textwidth,valign=c]{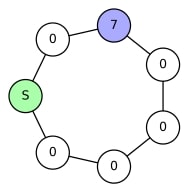}
        \caption{Example firing process for $C_n$ when $n=7$. The sink is colored green, vertices with 1 chip are colored blue, and vertices with $n=7$ chips are colored purple.}
        \label{figure: minimal compatibility example C7}
    \end{figure}
    
    Now we will look at even $n$.
    If $v_i$ is diametrically opposite the sink, we can send 2 chips from the sink to $v_i$, and as $n$ is even we can do this $\frac{n}{2}$ times to yield $v_i$ having $n$ chips and no other vertex having chips.
    Otherwise, from the sink, give each vertex except the vertex diametrically opposite the sink a single chip as in the odd case.
    Then, pair vertices equidistant from vertex $v_i$ except the vertex diametrically opposite $v_i$ and fire their chips to vertex $v_i$ like before.
    This results in vertex~$v_i$ having $n-1$ chips, the vertex diametrically opposite the sink having $-1$ chips, and the vertex diametrically opposite $v_i$ having 1 chip.
    Then, fire chips from the sink and the vertex with 1 chip to~$v_i$ and the vertex diametrically opposite the sink, resulting in $v_i$ having $n$ chips and no other vertex having chips.
    See \cref{figure: minimal compatibility example C6} for an example of the described firing process.
    
    \begin{figure}[h!]
        \centering
        \includegraphics[width=0.2\textwidth,valign=c]{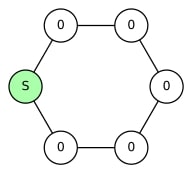}
        $\Rightarrow$
        \includegraphics[width=0.2\textwidth,valign=c]{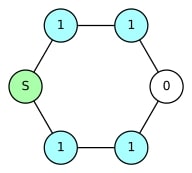}
        $\Rightarrow$
        \includegraphics[width=0.2\textwidth,valign=c]{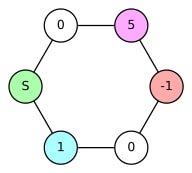}
        $\Rightarrow$
        \includegraphics[width=0.2\textwidth,valign=c]{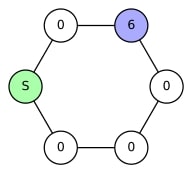}
        \caption{Example firing process for $C_n$ when $n=6$. The sink is colored green, vertices with 1 chip are colored blue, vertices with $-1$ chips are colored red, vertices with $n-1=5$ chips are colored pink, and vertices with $n=6$ chips are colored purple.}
        \label{figure: minimal compatibility example C6}
    \end{figure}
    
    Finally, we will show that $n$ is the smallest element in $I$.
    We will do this by showing the firing vector for the case where the sink and $v_i$ are adjacent is irreducible.
    Let the sink be vertex 1, with the vertices labeled in order so that $v_i$ is vertex $n$.
    First, fire vertex 1. Then, fire vertices 1 and 2. After that, fire vertices 1, 2, and 3, and so on, with the final step firing vertices 1 through $n-1$.
    Notice that each step pushes a chip in the positive direction; the first step sends a chip to vertex 2, the next step moves that chip to vertex 3, and so on. After all these steps, that chip will arrive at vertex $n$.
    At the same time, each step gives vertex $n$ a chip from the sink, vertex 1. Hence, after these $n-1$ steps, vertex $n$ will receive $1 + (n-1) = n$ chips, with no other vertex having chips.
    During this process, vertex $n-1$ was only fired once, and thus the firing vector is irreducible.
\end{proof}

\section{The Complete Identity Property}\label{section:CIP}
    
    Using the definition of compatibility, a graph has the complete maximal identity property if and only if $m_G$ is compatible.
    With this concept, the definition of the complete maximal identity property can be generalized to simply when there exists a compatible configuration that is recurrent for all sinks but need not be $m_G$, as seen in the following definition.
    
    \begin{definition}\label{def:CIP}
        A graph $G$ is said to have the \textbf{complete identity property} if there exists a chip configuration $c$ on all the vertices such that for all choices of sink $s$, the configuration $c$ with respect to sink $s$ is the recurrent identity for the sandpile of $G$ at $s$.
        
        
        If $G$ has the complete identity property, $c_G$ is the chip configuration on all vertices that gives the recurrent identity for all choices of sink $s$.
    \end{definition}
    
    Note that if a graph $G$ has the complete maximal identity property, it has the complete identity property.
    
    Odd cycle graphs have the complete maximal identity property, but even cycle graphs do not.
    The generalization of the complete maximal identity property to the complete identity property helps resolve this, as seen by \cref{proposition: CIP even cycle +}.
    
    \begin{proposition}\label{proposition: CIP even cycle +}
        For any positive integer $n$, attaching a single tree of size $1$ to any vertex in the even cycle graph $C_{2n}$ results in a graph with the complete identity property.
    \end{proposition}
    \begin{proof}
        Say $G$ is the graph resulting from attaching a single tree of size 1 to any vertex in $C_{2n}$.
        Let vertex $v_0$ be the vertex added to the even cycle graph and vertex $v_1$ be the vertex connected to it. 
        Let the vertices of the cycle graph be numbered~$v_1$ through $v_{2n}$ in clockwise order.
        Notice that $v_{n+1}$ is diametrically opposite of $v_1$.
        See \cref{subfigure: C6+} for an example of $G$, when $n=3$.
        
        We claim the common recurrent identity configuration $c_G$ is the configuration that has 0 chips at $v_0$ and $v_{n+1}$, 2 chips at $v_1$, and 1 chip everywhere else.
        Let this configuration be denoted as $d_G$ for the proof.
        We will show that this configuration is the recurrent identity for all of the vertices.
        See \cref{subfigure: dG C6+} for an example of $d_G$, when $n=3$.

        \begin{figure}[htbp!]
            \centering
            \begin{minipage}[t]{0.4\textwidth}
            \centering
            \includegraphics[width=0.75\textwidth]{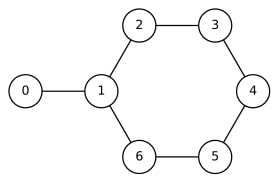}
            \caption{$G$ when $n=3$}
            \label{subfigure: C6+}
            \end{minipage}
            \quad\quad
            \begin{minipage}[t]{0.5\textwidth}
            \centering
            \includegraphics[width=0.6\textwidth]{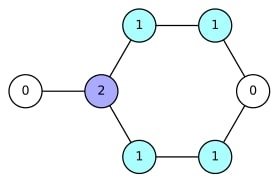}
            \caption{$d_G$ when $n = 3$. Vertices with 1 and 2 chips are colored blue and purple, respectively.}
            \label{subfigure: dG C6+}
            \end{minipage}
        \end{figure}
        
        \textbf{Case 1: sink at $v_0$.}
        To prove that the configuration $d_G$ is recurrent, take any stable configuration and add chips to it to yield the configuration that has 2 more chips than $d_G$ at $v_{n+1}$, which is the configuration resulting from taking the maximal stable configuration and adding 1 chip to~$v_{n+1}$, or $m_G + \mathbf{e}_{n+1}$.
        From here, we stabilize to reach $d_G$.
        Vertex $v_{n+1}$ fires, activating $v_{n+1\pm1}$.
        Vertex~$v_{n+1\pm1}$ fires, activating $v_{n+1}$ again as well as $v_{n+1\pm2}$.
        This process continues, sending the two chips originally at $v_{n+1}$ to $v_1$, which fires, activating the vertices all the way back to $v_{n+1}$, and then repeats this process again, giving the two chips to $v_0$, thus reaching $d_G$.
        Hence, $d_G$ is recurrent.
        Furthermore, using the same reasoning, all the chips in $d_G$ may be moved to $v_0$ as there are no chips at $v_{n+1}$, and the chips at $v_{n+1\pm i}$  (chips that are the same distance from vertex 1) may be paired together and moved towards $v_1$, and from there all the chips moved to $v_0$ by firing the entire cycle graph as many times as necessary.
        
        \textbf{Case 2: sink at $v_1$.}
        Follow the same procedure as $v_0$, ending when the two chips reach~$v_1$.
        
        \textbf{Case 3: sink at $v_{n+1}$.}
        Follow a similar procedure but for the proof of $d_G$ being recurrent, add chips to any stable configuration to yield the configuration that is the maximal stable configuration with two extra chips at $v_1$, and stabilize to yield $d_G$. 
        
        \textbf{Case 4: sink at $v_i$ for $2 \leq i \leq n$.}
        To prove that the configuration $d_G$ is recurrent, take any stable configuration and add chips to it to yield the configuration that has 1 more chip than $m_G$ at $v_{n+i}$.
        Firing $v_0$ immediately after each time $v_1$ fires, the proof of the recurrence of $d_G$ proceeds exactly like the case for the sink at $v_1$, acting on the cycle graph.
        
        To prove that $d_G$ is equivalent to the all-zero configuration, pair off vertices equidistant from~$v_i$ in the cycle graph and move their chips to $v_i$, resulting in the configuration with $v_1$ having 1 chip,~$v_{n+1}$ having $-1$ chips, and $v_{n+i}$ having 1 chip.
        Using this observation, we move the two chips at $v_1$ and $v_{n+i}$ away from each other until the chip from $v_{n+i}$ reaches $v_{n+1}$, and thus the chip at $v_1$ reaches $v_i$, the sink. This results in the all-zero configuration, and thus $d_G$ is the recurrent identity at~$v$.
        
        \textbf{Case 5: sink at $v_i$ for $n+2 \leq i \leq 2n$.}
        By symmetry, follow the same procedure as described in case 4 for $v_{2n + 2 - i}$ where $2 \leq 2n + 2 - i \leq n$.
        
        This completes the proof.
    \end{proof}
    
    In order to create a graph from an even cycle graph that has the complete identity property, we needed to add a single tree of size 1 to any vertex.
    This property is also seen in the case of complete bipartite graphs, as the following theorem shows.
    
    \begin{theorem}\label{theorem: CIP K_mn +}
        For all positive integers $m,n$, attaching a single tree of size $1$ to any vertex in the complete bipartite graph $K_{m,n}$ results in a graph that has the complete identity property.
    \end{theorem}
    \begin{proof}
        
        It suffices to show the result for when the additional tree is attached to one of the $m$ vertices, as the same proof would hold with $m$ and $n$ interchanged. Furthermore, we may assume that $m,n\geq 2$, as if any of them are 1, the resulting graph is a tree and thus has the complete maximal identity property and hence the complete identity property by \cref{proposition: CMIP complete odd cycle tree}.
        
        Say the additional vertex is $a$, the vertex it is attached to is $b$, the $n$ vertices of one side of the complete bipartite graph compose set $C$, and the $m-1$ vertices in the set of $m$ vertices of one side of the complete bipartite graph that is not $b$ compose set $D$.
        We claim $c_G$ is the configuration with $a$ having 0 chips, $b$ having $n$ chips, each vertex in $C$ having $m-1$ chips, and each vertex in $D$ having 0 chips.
        Let this configuration be denoted as $d_G$ for the proof. See \cref{fig:K5_7+} for an example of the graph and the naming convention for the vertices.
        
        \begin{figure}[h!]
        \centering

        \begin{tikzpicture}[thick,scale=0.3,font=\Large,baseline]
            \foreach \x in {0,...,6} {
                \foreach \y in {0,...,4} {
                    \draw (\x*4,0) -- (\y*6,10);
                }
            }
            \draw (0,10) -- (-5,5);
            \foreach \x in {0,...,6}
                \filldraw[color=black,fill=white] (\x*4,0) circle (40pt);
            \foreach \x in {0,...,4}
                \filldraw[color=black,fill=white] (\x*6,10) circle (40pt);
            \filldraw[color=black,fill=white] (-5,5) circle (40pt);
            \foreach \x in {5,...,11}
                \node at (\x*4-20,0) {$v_{\x}$};
            \foreach \x in {0,...,4}
                \node at (\x*6,10) {$v_{\x}$};
            \node at (-5,5) {$v_{12}$};
        \end{tikzpicture}
        \caption{Example of $K_{m,n}$ with a tree of size 1 attached, for $m=5,n=7$. Vertex 12 is $a$, vertex 0 is $b$, vertices 5 through 11 compose set $C$, and vertices 1 through 4 compose set $D$.}\label{fig:K5_7+}
        
        \end{figure}
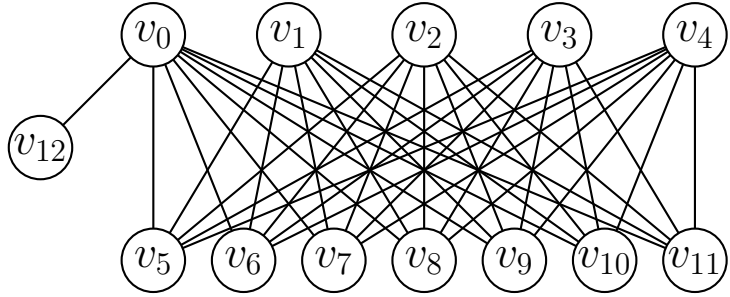
        
        We will separate into four cases depending on where the sink is.
        
        \textbf{Case 1: sink at $a$.}
        From any stable configuration add chips to reach the configuration with~$b$ having $n$ chips, each vertex in $C$ having $m-1$ chips, and each vertex in $D$ having $n$ chips.
        We then stabilize.
        This results in all vertices in $D$ firing, then all vertices in $C$ firing, then all vertices in $D$ firing again to result in $b$ having $2n$ chips and $C$ having $2m-3$ chips each.
        Then, firing $b$, all vertices in~$C$, and then all vertices in $D$ results in a net loss of 1 chip at $b$, and this occurs as long as $b$ has at least~$n+1$ chips to start the cycle, and the cycle starts with all vertices in $C$ having at least~$m-1$ chips, which it does as $m \geq 2$, so $2m-3 \geq m-1$.
        Hence, after $n$ iterations of this cycle, we eventually reach $b$ having $n$ chips, all vertices in $C$ having $2m-3$ chips each, and all vertices in~$D$ having 0 chips each.
        Fire all vertices in $C$ and then all vertices in $D$ to yield all vertices in $C$ having lost 1 chip each, and $b$ having gained $n$ chips.
        Follow the cycle previously shown to return~$b$ back to having $n$ chips.
        This process can thus result in a net loss of 1 chip at each vertex in $C$, as long as all vertices in $C$ had at least $m$ chips to start with, giving it at least $m-1$ chips when it enters the other cycle.
        Cycling this process until it is no longer possible results in $b$ having $n$ chips, all vertices in $C$ having $m-1$ chips each, and all vertices in $D$ having no chips.
        This is $d_G$, and it is stable.
        Hence, $d_G$ is recurrent.
        Now, to prove that $d_G$ is the recurrent identity, fire all vertices in~$C\cup D$ each $m-1$ times to clear all vertices in $C \cup D$ of chips, and then fire $b$, all vertices in $C$, and all vertices in $D$ until $b$ is also clear of chips.
        \cref{table:Kmn-a} shows a stabilization of this configuration to $d_G$, and then the unrestricted firing sequence that takes $d_G$ to $\mathbf{0}$, tracking the chips at each non-sink vertex.
        As during the process the vertices in $C$ are indistinguishable, the column for $C$ tracks the number of chips at each of the vertices in $C$.
        The chips at each vertex in $D$ are similarly tracked.
        Hence, $d_G$ is the recurrent identity.
        \begin{table}[htbp!]
        \centering
        \begin{tabular}{| c | c | c | c || c |}
            \hline
            $a$ &  $b$ & $C$ & $D$ & Firing step \\
            \hline
            Sink & $n$ & $m-1$ & $n$ & $D,C,D$ fire \\
            Sink & $2n$ & $2m-3$ & 0 & $b, C, D$ fire \\
            Sink & $2n-1$ & $2m-3$ & 0 & Repeat firing $b,C,D$ a total of $n$ times \\
            Sink & $n$ & $2m-3$ & 0 & $C,D$ fire, $b,C,D$ fire $n$ times. Do this $m-2$ times \\
            Sink & $n$ & $m-1$ & 0 & Reached $d_G$ \\
            \hline
            \hline
            Sink & $n$ & $m-1$ & 0 & $C,D$ fire $m-1$ times \\
            Sink & $nm$ & 0 & 0 & $b,C,D$ fire $nm$ times \\
            Sink & 0 & 0 & 0 & Reached $\mathbf{0}$ \\
            \hline

        \end{tabular}
        \caption{Firing sequences for proof that $d_G$ is the recurrent identity with sink at $a$}
        \label{table:Kmn-a}
        \end{table}       
        
        \textbf{Case 2: sink at $b$.}
        From any stable configuration add chips to reach the configuration with~$a$ having 0 chips, each vertex in $C$ having $m-1$ chips, and each vertex in $D$ having $n$ chips.
        Stabilizing this configuration results in first firing all vertices in $D$, resulting in all vertices in $C$ having $2m-2$ chips each.
        Then, firing all vertices in $C$ and then all vertices in $D$ results in a net loss of 1 chip at each vertex in $C$, and the cycle works as long as all vertices in $C$ start with at least $m$ chips so that they can fire. Repeating this process, we find that the stabilization of the specified accessible configuration is $d_G$. Hence, $d_G$ is recurrent. Now, to prove that $d_G$ is the recurrent identity, fire all vertices in $C \cup D$ each $m-1$ times to result in the all-zero configuration.
        \cref{table:Kmn-b} shows the firing sequences that prove $d_G$ is the recurrent identity for sink at $b$.
        \begin{table}[htbp!]
        \centering
        \begin{tabular}{| c | c | c | c || c |}
            \hline
            $a$ & $b$ & $C$ & $D$ & Firing step \\
            \hline
            0 & Sink & $m-1$ & $n$ & $D$ fires \\
            0 & Sink & $2m-2$ & 0 & $C, D$ fire repeatedly until no longer possible\\
            0 & Sink & $m-1$ & 0 & Reached $d_G$ \\
            \hline
            \hline
            0 & Sink & $m-1$ & 0 & $C,D$ fire $m-1$ times \\
            0 & Sink & 0 & 0 & Reached $\mathbf{0}$ \\
            \hline
        \end{tabular}
        \caption{Firing sequences for proof that $d_G$ is the recurrent identity with sink at $b$}
        \label{table:Kmn-b}
        \end{table}
        
        \textbf{Case 3: sink in $C$.}
        Say the sink is $s_C \in C$, and let $C' = C \setminus \{s_C\}$.
        From any stable configuration add chips to reach the configuration with $a$ having 0 chips, $b$ having $n$ chips, each vertex in $C'$ having $m$ chips, and each vertex in $D$ having $n$ chips.
        We fire all vertices in $C'$ to yield~$b$ and each vertex in $D$ having~$2n-1$ chips each, both being active.
        Notice that we may fire~$b$, $a$, and all vertices in $D$, and then all vertices in $C'$, as long as $b$ and all vertices in $D$ have at least~$n+1$ chips each, and without any requirement of the starting chips for all vertices in $C'$, as they will each get the $m$ chips they need to fire.
        This process results in a net loss of 1 chip at each of $b$ and all vertices in $D$.
        Performing this operation until it is no longer possible, with $n$ chips at~$b$ and each vertex in $D$ each, we then fire all vertices in $D$ to yield $d_G$. Hence, $d_G$ is recurrent.
        To prove that $d_G$ is the recurrent identity, notice that~$d_G$ is equivalent to the configuration that preceded it in the proof of its recurrence, which had $n$ chips at each of $b$ and all vertices in $D$; hence, backfiring $D$ results in this configuration.
        We now backfire $s_C$ a total of $n$ times to clear $b$ and all vertices in $D$ of chips, resulting in the all-zero configuration.
        \cref{table:Kmn-C} shows the firing sequences that prove $d_G$ is the recurrent identity for sink in $C$.
        \begin{table}[htbp!]
        \centering
        \begin{tabular}{| c | c | c | c || c |}
            \hline
            $a$ & $b$ & $C'$ & $D$ & Firing step \\
            \hline
            0 & $n$ & $m$ & $n$ & $C'$ fires \\
            0 & $2n-1$ & 0 & $2n-1$ & $b,a,D,C'$ fire repeatedly until no longer possible\\
            0 & $n$ & 0 & $n$ & $D$ fires \\
            0 & $n$ & $m-1$ & 0 & Reached $d_G$ \\
            \hline
            \hline
            0 & $n$ & $m-1$ & 0 & Backfire $D$ \\
            0 & $n$ & 0 & $n$ & $s_C$ backfires $n$ times \\
            0 & 0 & 0 & 0 & Reached $\mathbf{0}$ \\
            \hline
        \end{tabular}
        \caption{Firing sequence for proof of recurrence of $d_G$ with sink in $C$}
        \label{table:Kmn-C}
        \end{table}

        \textbf{Case 4: sink in $D$.}
        Say the sink is $s_D \in D$, and let $D' = D \setminus \{s_D\}$.
        From any stable configuration add chips to reach the configuration with $a$ having 0 chips, $b$ having $2n$ chips, each vertex in $C$ having $m-1$ chips, and each vertex in~$D'$ having $n$ chips.
        Notice that at this state $b$ and all vertices in~$D'$ are active.
        Furthermore, firing these vertices and then $a$ gives $m-1$ chips to each vertex in $C$, so as long as they each have at least 1 chip, they can fire.
        This results in $b$ and~$D'$ regaining their $n$ chips, and costing each vertex in $C$ $m$ chips, meaning that each vertex in~$C$ has lost one chip in total.
        We perform this procedure $m-1$ times to result in the configuration with~$2n$ chips at $b$ and $n$ chips at each vertex in $D'$.
        Firing $b$, $a$, and all vertices in $D'$ then gives~$d_G$. Hence,~$d_G$ is recurrent.
        To prove that $d_G$ is the recurrent identity, fire $a$ and $b$ to clear $a$ and~$b$ of chips, resulting in $C$ having $m$ chips. Backfire the sink $m$ times to result in the all-zero configuration.
        \cref{table:Kmn-D} shows the firing sequences that prove $d_G$ is the recurrent identity for sink in $D$.
        \begin{table}[htbp!]
            \centering
            \begin{tabular}{| c | c | c | c || c |}
            \hline
            $a$ & $b$ & $C$ & $D'$ & Firing step \\
            \hline
            0 & $2n$ & $m-1$ & $n$ & $b,D',a,C$ fire $m-1$ times. \\
            0 & $2n$ & 0 & $n$ & $b,D',a$ fire \\
            0 & $n$ & $m-1$ & 0 & Reached $d_G$ \\
            \hline
            \hline
            0 & $n$ & $m-1$ & 0 & $a,b$ fire \\
            0 & 0 & $m$ & 0 & $s_D$ backfires $m$ times \\
            0 & 0 & 0 & 0 & Reached $\mathbf{0}$ \\
            \hline
        \end{tabular}
        \caption{Firing sequence for proof of recurrence of $d_G$ with sink in $D$}
        \label{table:Kmn-D}
        \end{table}
        
        This completes the proof.
    \end{proof}
    
    By \cref{proposition: CIP even cycle +} and \cref{theorem: CIP K_mn +}, we see that even cycles and complete bipartite graphs both have the property that attaching a single tree of size 1, or a single edge and vertex, to any vertex in the graph results in a graph with the complete identity property.
    In addition, adding a single edge and vertex to a tree results in a tree, which has the complete maximal identity property by \cref{proposition: CMIP complete odd cycle tree}, which implies the complete identity property.
    All of these three types of graphs are bipartite.
    While not all bipartite graphs, or even regular bipartite graphs, have this property (for example, the hypercube in 3 dimensions, which is isomorphic to $C_4 \cp K_2$, where the Cartesian product $\cp$ is defined in \cref{section:Conjectures}, is one such counterexample), a computer search found that all connected graphs of 10 vertices or less which had this property were bipartite.
    This motivates the following conjecture.
    
    \begin{conjecture}\label{conjecture: CIP+ >>> bipartite}
        Let $G$ be a connected graph. If for all vertices $v\in G$, attaching a single tree of size $1$ to $v$ results in a graph with the complete identity property, then $G$ is bipartite.
    \end{conjecture}
    

\section{Conjectures on Graph Products}\label{section:Conjectures}
\begin{figure}[h!]
        \centering
        \begin{subfigure}[b]{0.2\textwidth}
            \centering
            \begin{tikzpicture}[thick,baseline]
                \draw (0,0) -- (2,0);
                \foreach \x in {0,1,2}
                    \filldraw[color=black,fill=white] (\x,0) circle (6pt);
            \end{tikzpicture}
            \caption{$P_3$}
            \label{subfigure:P3}
        \end{subfigure}
        \quad
        \begin{subfigure}[b]{0.2\textwidth}
            \centering
            \begin{tikzpicture}[thick,baseline]
                \foreach \x in {0,1,2} {
                    \draw (\x,0) -- (\x,2);
                    \draw (0,\x) -- (2,\x);
                }
                \foreach \x in {0,1,2} {
                    \foreach \y in {0,1,2}
                        \filldraw[color=black,fill=white] (\x,\y) circle (6pt);
                }
            \end{tikzpicture}
            \caption{$P_3 \protect\cp P_3$}
            \label{subfigure:P3cpP3}
        \end{subfigure}
        \quad
        \begin{subfigure}[b]{0.2\textwidth}
            \centering
            \begin{tikzpicture}[thick,baseline]
                \draw (0,1) -- (1,0);
                \draw (0,2) -- (2,0);
                \draw (1,2) -- (2,1);
                \draw (0,1) -- (1,2);
                \draw (0,0) -- (2,2);
                \draw (1,0) -- (2,1);
                \foreach \x in {0,1,2} {
                    \foreach \y in {0,1,2}
                        \filldraw[color=black,fill=white] (\x,\y) circle (6pt);
                }
            \end{tikzpicture}
            \caption{$P_3 \times P_3$}
            \label{subfigure:P3tpP3}
        \end{subfigure}
        \quad
        \begin{subfigure}[b]{0.2\textwidth}
            \centering
            \begin{tikzpicture}[thick,baseline]
                \draw (0,1) -- (1,0);
                \draw (0,2) -- (2,0);
                \draw (1,2) -- (2,1);
                \draw (0,1) -- (1,2);
                \draw (0,0) -- (2,2);
                \draw (1,0) -- (2,1);
                \foreach \x in {0,1,2} {
                    \draw (\x,0) -- (\x,2);
                    \draw (0,\x) -- (2,\x);
                }
                \foreach \x in {0,1,2} {
                    \foreach \y in {0,1,2}
                        \filldraw[color=black,fill=white] (\x,\y) circle (6pt);
                }
            \end{tikzpicture}
            \caption{$P_3 \protect\sp P_3$}
            \label{subfigure:P3spP3}
        \end{subfigure}
        \caption{The Cartesian, tensor, and strong products between two copies of $P_3$}
        \label{figure: cp tp sp}
    \end{figure}
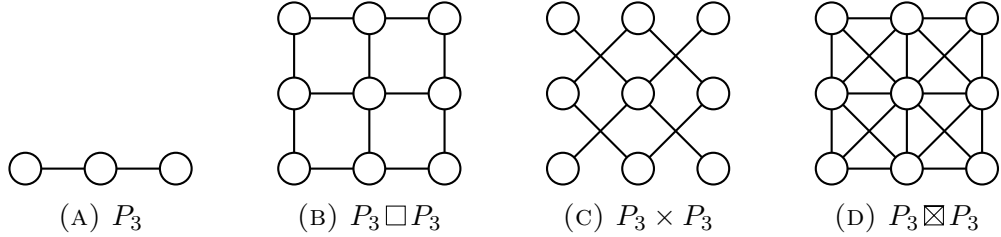

    Recall that the Cartesian product, tensor product, and strong product are binary operations on graphs that form a graph whose vertices are ordered pairs of vertices of the two daughter graphs.
    For the Cartesian product, denoted $\cp$, two vertices share an edge if in one of the daughter graphs the two vertices share an edge and in the other the vertices are the same.
    For the tensor product, denoted $\times$, two vertices share an edge if in both daughter graphs the two vertices share an edge.
    The strong product, denoted $\sp$, is the union of the Cartesian product and the tensor product.
    See \cref{figure: cp tp sp} for examples of the Cartesian, tensor, and strong products between two copies of the path graph with 3 vertices $P_3$.
    
    Investigating the behavior of graph products with respect to the complete maximal identity property, we find that in general, the strong product, Cartesian product, and tensor product all do not preserve the complete maximal identity property.
    However, we do find some patterns.
    
    \begin{proposition}\label{proposition:S-PxP}
        Let $P_n$ be the path graph with $n$ vertices. The strong product between $P_2$ and $P_k$ has the complete maximal identity property if and only if $k=2$ or $k \equiv 1 \pmod{3}$.
    \end{proposition}
    \begin{proof}
        
        We will refer to the canonical labelling of the vertices of the graph $P_2 \sp P_k$ as via the ordered pairs $(i,j)$ where $i\in\{0,1\}$ and $j\in\{0,1,\dots,k-1\}$. By symmetry, it suffices to prove the graph has the maximal identity property for sinks with $i=0$ and $j\leq \frac{k-1}{2}$.
        

        We will prove the $j=0$ case first.
        Starting with the maximal stable configuration, incrementally fire all the vertices with second coordinate~$k-1$ until those have no chips, and then fire all the vertices with second coordinate at least~$k-2$ until those have no chips, and so on.
        The resulting configuration has~$4k-6$ chips at each of the two vertices with second coordinate~$1$.
        By symmetry, in order for the two vertices with second coordinate 1 to have the same number of chips (eventually~0), the number of times the vertices with second coordinate at least 1 can be fired must all be the same.
        Fire all of these vertices $k-2$ times.
        Then the three vertices connected to the sink each have $2k-2$ chips.
        Fire all non-sink vertices $2k-2$ times to clear the graph of all chips.
        
        Now we assume $k\geq3$.
        For $1<j\leq \frac{k-1}{2}$, we follow a similar process to result in the vertices with second coordinate~$j-1$ having all the chips originally with second coordinate less than $j$, and thus having $4j-2$ chips each.
        Similarly, the vertices with second coordinate $j+1$ have $4k-4j-6$ chips each. Fire the vertices with second coordinate at least $j+1$ a total of $k-2j-1$ times each to have the vertices with second coordinate $j\pm1$ having $4j-2$ chips each, and $(1,j)$ having $4k-8j$ chips.
        
        In order for the vertices with second coordinate not equal to $j$ to all have 0 chips, they must all be fired the same number of times.
        Hence, we essentially only have two operations: fire $(1,j)$ or fire all vertices with second coordinate not equal to $j$.
        From this, we can see that we can fire all the non-sink vertices $4j-2$ times, resulting in all the non-sink vertices having no chips except~$(1,j)$ having $4k-12j+2$ chips.
        Now notice that the smallest increment by which we can change the number of chips at that vertex without changing the other vertices is to fire it twice and then fire all the non-sink vertices once.
        Or, in other words, fire it once and backfire the sink once. This results in a net loss of 6 chips.
        So for $P_2 \sp P_k$ to have the complete maximal identity property we must have that $4k-12j+2\equiv 0 \mod{6}$ for all values of $j$, or equivalently $k \equiv 1 \mod{3}$.
        If this were not the case, we would be able to get the all-zero configuration with a fractional number of firings, and as the reduced Laplacian is nonsingular, this is the unique firing vector needed, and thus the maximal stable configuration is not equivalent to the all-zero configuration.
        
        To prove that the other graphs do not have the complete maximal identity property, let $k \not\equiv 1 \mod{3}$, and use the same process as before to arrive at $(1,j)$ having 2 or 4 chips.
        Using the fact that the reduced Laplacian is non-singular, there is a unique firing vector that results in this configuration.
        However, it does not have integer entries, as the configuration with 6 chips at $(1,j)$ does not have all of its entries being multiples of 3, rather having all of the vertices with second coordinate not equal to $j$ being backfired once and $(1,j)$ being backfired twice.
        Hence, the two configurations resulting from $k \not\equiv 1 \mod{3}$ are not equivalent to the identity, and thus do not have the complete maximal identity property.
    \end{proof}
    
    After looking at whether $P_i \sp P_j$ has the complete maximal identity property for all values of $i$ and $j$ where~$1 < i,j \leq 100$, we conjecture that the cases presented in \cref{proposition:S-PxP} are the only such graphs with the complete maximal identity property:
    
    \begin{conjecture}\label{conjecture:S-PxP}
        For $i \leq j \in \mathbb{Z}_{>1}$, the only graphs $P_i \sp P_j$ which have the complete maximal identity property are $P_2 \sp P_j$ where $j\equiv 1 \mod{3}$ or $j=2$, which yields $K_4$.
    \end{conjecture}
    
    A computer program in SageMath verified that the conjecture holds for $1 < i \leq j \leq 100$. 
    
    Similarly, the following proposition on the Cartesian product was proven.
    
    \begin{proposition}\label{proposition:C-KxP}
        The Cartesian product between $K_4$ and $P_2$ has the complete maximal identity property.
    \end{proposition}
    
    After looking at whether or not $K_i \cp P_j$ has the complete maximal identity property for all values of $i$ and $j$ where~$1 < i,j \leq 50$, we conjecture that that the case presented in \cref{proposition:C-KxP} is the only such graph with the complete maximal identity property.
    
    \begin{conjecture}\label{conjecture:C-KxP}
        For $i,j \in \mathbb{Z}_{>1}$, the only graph $K_i \cp P_j$ which has the complete maximal identity property is $K_4 \cp P_2$.
    \end{conjecture}
    
    Another computer program in SageMath verified that the conjecture holds for $1 < i,j \leq 50$.

    We provide a proof of the special cases of \cref{conjecture:C-KxP} when $j=2,3$.
    
    \begin{proof}[Proof of the $j=2,3$ cases]


        For $j=2$, let the sink be (0,0).
        Instead of looking at abscissas of 0 through $i-1$, by symmetry all the positive abscissas for a specified ordinate must fire the same number of times (even if this were not true, the solution that follows would yield a non-integer number of firings for each vertex, and as the reduced Laplacian is non-singular, no other firing vector will yield the all-zero configuration, and thus the proof holds).
        So, we will combine the vertices with the same ordinate and a positive abscissa together to yield a weighted cycle graph of 4 vertices. The reduced Laplacian is
        \begin{equation*}
            \Delta' = \begin{bmatrix}
            i & 1-i & 0 \\
            1-i & 2i-2 & 1-i \\
            0 & 1-i & 2i-2 \\
            \end{bmatrix}.
        \end{equation*}
        Thus, the firing vector $\mathbf{v}$ to reach the maximal stable configuration $\mathbf{c}$ is
        \begin{equation*}
            \mathbf{v} =(\Delta')^{-1} \mathbf{c} = (\Delta')^{-1} \begin{bmatrix} i - 1 \\ (i-1)^2 \\ (i-1)^2 \end{bmatrix} = \frac{1}{i+2}\begin{bmatrix} 3i(i-1) \\ (i-1)(3i+2) \\ 2(i-1)(i+1) \end{bmatrix}.
        \end{equation*}
        To have the complete maximal identity property, it suffices to show that for this particular sink, $\mathbf{v}$ has integer entries. So
        \[ i+2 \mid \gcd\left( 3i(i-1), (i-1)(3i+2), 2(i-1)(i+1)\right). \]
        
        We will first analyze $i+2 \mid 3i(i-1)$.
        Notice that $\gcd(i,i+2)=\gcd(2,i)\mid 2$. 
        We also have~$\gcd(i+2,i-1)=\nolinebreak\gcd(i+2,3)\mid 3$ and $\gcd(i+2,3)\mid 3$.
        Hence, we find that $i+2 \mid 2\cdot 3^2$.
        This corresponds to integer values of $i$ greater than 1 being 4, 7, and 16.
        Verifying that these hold for the other two divisibility criteria, we find only $i=4$ yields the complete maximal identity property.
        
        For $j = 3$, we follow the same process for combining vertices. The reduced Laplacian is 
        \begin{equation*}
            \Delta' = \begin{bmatrix}
            i   & 1-i  & -1   & 0    & 0\\
            1-i & 2i-2 & 0    & 1-i  & 0\\
            -1  & 0    & i+1  & 1-i  & 0\\
            0   & 1-i  & 1-i  & 3i-3 & 1-i\\
            0   & 0    & 0    & 1-i  & 2i-2\\
            \end{bmatrix}.
        \end{equation*}
        The firing vector $\mathbf{v}$ to reach the maximal stable configuration $c$ is
        \begin{equation*}
            \mathbf{v}=(\Delta')^{-1}\mathbf{c} = \frac{1}{(i+1)(i+3)}
            \begin{bmatrix}
            2i(3i-2)(i+3) \\
            (3i-2)(2i^2 +6i + 1) \\
            5i^3 + 8i^2 - 6i + 1 \\
            5i^3 + 11i^2 -5i -1 \\
            (3i-2)(i^2+3i+1) \\
            \end{bmatrix}.
        \end{equation*}
        This requires that
        \[ (i+1)(i+3)\mid \gcd\left( 5i^3 + 8i^2 - 6i + 1, 5i^3 + 11i^2 -5i -1 \right)\]
        which results in
        \[ (i+1)(i+3)\mid 3i^2 + i -2. \]
        With $i+1\mid 3i^2+i-2$, the condition is equivalent to $i+3 \mid 3i-2$, or $i+3 \mid 11$, or $i=8$.
        But this does not yield a vector with integer entries, so hence there are no solutions for $j=3$.
    \end{proof}


\section*{Acknowledgements}
    We would like to thank Professor David Perkinson of the Department of Mathematics at Reed College for suggesting the project and his guidance during the chip-firing research project. We would also like to thank Dr. Tanya Khovanova of the Department of Mathematics at the Massachusetts Institute of Technology (MIT) for her advice throughout the research process. We thank Dr. Slava Gerovitch and Professor Pavel Etingof of the Department of Mathematics at MIT for their operation of the PRIMES-USA program; we also thank the MIT PRIMES-USA Program as well as the Department of Mathematics at MIT.


\bibliographystyle{plain}
\bibliography{ref}

\begin{thebibliography}{10}

\bibitem{bak2013nature}
Per Bak.
\newblock {\em How nature works: the science of self-organized criticality}.
\newblock Springer Science \& Business Media, 2013.

\bibitem{bak1990forest}
Per Bak, Kan Chen, and Chao Tang.
\newblock A forest-fire model and some thoughts on turbulence.
\newblock {\em Physics letters A}, 147(5-6):297--300, 1990.

\bibitem{bak1987self}
Per Bak, Chao Tang, and Kurt Wiesenfeld.
\newblock Self-organized criticality: An explanation of 1/f noise.
\newblock {\em Phys. Rev. Lett.}, 59:381--384, 1987.

\bibitem{bartolozzi2006scale}
Marco Bartolozzi, Derek~B Leinweber, and Anthony~W Thomas.
\newblock Scale-free avalanche dynamics in the stock market.
\newblock {\em Physica A: Statistical Mechanics and its Applications},
  370(1):132--139, 2006.

\bibitem{biggs1999chip}
N.~L. Biggs.
\newblock Chip-firing and the critical group of a graph.
\newblock {\em J. Algebraic Combin.}, 9(1):25--45, 1999.

\bibitem{biggs1997algebraic}
Norman~L. Biggs.
\newblock Algebraic potential theory on graphs.
\newblock {\em Bull. London Math. Soc.}, 29(6):641--682, 1997.

\bibitem{biondo2015modeling}
Alessio~Emanuele Biondo, Alessandro Pluchino, and Andrea Rapisarda.
\newblock Modeling financial markets by self-organized criticality.
\newblock {\em Physical Review E}, 92(4):042814, 2015.

\bibitem{bjorner1992chip}
Anders Bj\"{o}rner and L\'{a}szl\'{o} Lov\'{a}sz.
\newblock Chip-firing games on directed graphs.
\newblock {\em J. Algebraic Combin.}, 1(4):305--328, 1992.

\bibitem{bjorner1991chip}
Anders Bj\"{o}rner, L\'{a}szl\'{o} Lov\'{a}sz, and Peter~W. Shor.
\newblock Chip-firing games on graphs.
\newblock {\em European J. Combin.}, 12(4):283--291, 1991.

\bibitem{burridge1967model}
Robert Burridge and Leon Knopoff.
\newblock Model and theoretical seismicity.
\newblock {\em Bulletin of the seismological society of america},
  57(3):341--371, 1967.

\bibitem{caracciolo2008explicit}
Sergio Caracciolo, Guglielmo Paoletti, and Andrea Sportiello.
\newblock Explicit characterization of the identity configuration in an abelian
  sandpile model.
\newblock {\em J. Phys. A}, 41(49):495003, 17, 2008.

\bibitem{caracciolo2012multiple}
Sergio Caracciolo, Guglielmo Paoletti, and Andrea Sportiello.
\newblock Multiple and inverse topplings in the abelian sandpile model.
\newblock {\em The European Physical Journal Special Topics}, 212(1):23--44,
  2012.

\bibitem{chialvo2010emergent}
Dante~R Chialvo.
\newblock Emergent complex neural dynamics.
\newblock {\em Nature physics}, 6(10):744--750, 2010.

\bibitem{cori2000sandpile}
Robert Cori and Dominique Rossin.
\newblock On the sandpile group of dual graphs.
\newblock {\em European J. Combin.}, 21(4):447--459, 2000.

\bibitem{corry2018divisors}
Scott Corry and David Perkinson.
\newblock {\em Divisors and sandpiles}.
\newblock American Mathematical Society, Providence, RI, 2018.

\bibitem{creutz1991abelian}
Michael Creutz.
\newblock Abelian sandpiles.
\newblock {\em Computers in Physics}, 5(2):198--203, 1991.

\bibitem{dhar1990self}
Deepak Dhar.
\newblock Self-organized critical state of sandpile automation models.
\newblock {\em Phys. Rev. Lett.}, 64(23):2837, 1990.

\bibitem{dhar1999abelian}
Deepak Dhar.
\newblock The abelian sandpile and related models.
\newblock {\em Physica A: Statistical Mechanics and its Applications},
  263(1-4):4--25, 1999.

\bibitem{dhar2006theoretical}
Deepak Dhar.
\newblock Theoretical studies of self-organized criticality.
\newblock {\em Physica A: Statistical Mechanics and its Applications},
  369(1):29--70, 2006.

\bibitem{dhar1994inverse}
Deepak Dhar and SS~Manna.
\newblock Inverse avalanches in the abelian sandpile model.
\newblock {\em Physical Review E}, 49(4):2684, 1994.

\bibitem{dhar1995algebraic}
Deepak Dhar, Philippe Ruelle, Siddhartha Sen, and D-N Verma.
\newblock Algebraic aspects of abelian sandpile models.
\newblock {\em Journal of physics A: mathematical and general}, 28(4):805,
  1995.

\bibitem{engel1975probabilistic}
Arthur Engel.
\newblock The probabilistic abacus.
\newblock {\em Educational studies in mathematics}, 6(1):1--22, 1975.

\bibitem{engel1976does}
Arthur Engel.
\newblock Why does the probabilistic abacus work?
\newblock {\em Educational Studies in Mathematics}, 7(1-2):59--69, 1976.

\bibitem{hesse2014self}
Janina Hesse and Thilo Gross.
\newblock Self-organized criticality as a fundamental property of neural
  systems.
\newblock {\em Frontiers in systems neuroscience}, 8:166, 2014.

\bibitem{holroyd2008chip}
Alexander~E. Holroyd, Lionel Levine, Karola M\'{e}sz\'{a}ros, Yuval Peres,
  James Propp, and David~B. Wilson.
\newblock Chip-firing and rotor-routing on directed graphs.
\newblock In {\em In and out of equilibrium. 2}, volume~60 of {\em Progr.
  Probab.}, pages 331--364. Birkh\"{a}user, Basel, 2008.

\bibitem{klivans2018mathematics}
Caroline~J Klivans.
\newblock {\em The mathematics of chip-firing}.
\newblock Chapman and Hall/CRC, 2018.

\bibitem{le2002identity}
Yvan Le~Borgne and Dominique Rossin.
\newblock On the identity of the sandpile group.
\newblock {\em Discrete Math.}, 256(3):775--790, 2002.
\newblock LaCIM 2000 Conference on Combinatorics, Computer Science and
  Applications (Montreal, QC).

\bibitem{levine2009sandpile}
Lionel Levine.
\newblock The sandpile group of a tree.
\newblock {\em European J. Combin.}, 30(4):1026--1035, 2009.

\bibitem{newman1996self}
MEJ Newman.
\newblock Self-organized criticality, evolution and the fossil extinction
  record.
\newblock {\em Proceedings of the Royal Society of London. Series B: Biological
  Sciences}, 263(1376):1605--1610, 1996.

\bibitem{pu2013developing}
Jiangbo Pu, Hui Gong, Xiangning Li, and Qingming Luo.
\newblock Developing neuronal networks: self-organized criticality predicts the
  future.
\newblock {\em Scientific reports}, 3(1):1--6, 2013.

\bibitem{scheinkman1994self}
Jose~A Scheinkman and Michael Woodford.
\newblock Self-organized criticality and economic fluctuations.
\newblock {\em The American Economic Review}, 84(2):417--421, 1994.

\bibitem{sneppen1995evolution}
Kim Sneppen, Per Bak, Henrik Flyvbjerg, and Mogens~H Jensen.
\newblock Evolution as a self-organized critical phenomenon.
\newblock {\em Proceedings of the National Academy of Sciences},
  92(11):5209--5213, 1995.

\bibitem{sole1996extinction}
Ricard~V Sol{\'e} and Susanna~C Manrubia.
\newblock Extinction and self-organized criticality in a model of large-scale
  evolution.
\newblock {\em Physical Review E}, 54(1):R42, 1996.

\bibitem{stanley2013algebraic}
Richard~P. Stanley.
\newblock {\em Algebraic combinatorics}.
\newblock Undergraduate Texts in Mathematics. Springer, Cham, 2018.

\bibitem{tardos1988polynomial}
G\'{a}bor Tardos.
\newblock Polynomial bound for a chip firing game on graphs.
\newblock {\em SIAM J. Discrete Math.}, 1(3):397--398, 1988.

\bibitem{turcotte2004landslides}
Donald~L Turcotte and Bruce~D Malamud.
\newblock Landslides, forest fires, and earthquakes: examples of self-organized
  critical behavior.
\newblock {\em Physica A: Statistical Mechanics and its Applications},
  340(4):580--589, 2004.

\bibitem{van2001algorithmic}
Jan van~den Heuvel.
\newblock Algorithmic aspects of a chip-firing game.
\newblock {\em Combin. Probab. Comput.}, 10(6):505--529, 2001.

\bibitem{wagner2000critical}
David~G Wagner.
\newblock The critical group of a directed graph.
\newblock {\em arXiv preprint math/0010241}, 2000.

\end{thebibliography}


\end{document}